\documentclass[12pt]{amsart}
\usepackage{latexsym,amssymb,amsmath}
\usepackage{amsmath,amsfonts}
\usepackage{epsfig}
\usepackage{graphicx}
\usepackage{verbatim}
\usepackage{tikz}
\usepackage{enumitem}
\usetikzlibrary{calc,decorations.markings}

\newtheorem{theorem}{Theorem}[section]
\newtheorem{lemma}[theorem]{Lemma}

\newtheorem{definition}{Definition}
\allowdisplaybreaks

\newcommand{\N}{\mathbb{N}}
\newcommand{\R}{\mathbb{R}}

\newcommand{\T}{\mathbb{T}}

\newcommand{\lb}{\langle}
\newcommand{\rb}{\rangle}
\newcommand{\I}{\operatorname{I}}
\newcommand{\II}{\operatorname{II}}
\newcommand{\III}{\operatorname{III}}
\renewcommand{\d}{\operatorname{d}\!}

\renewcommand{\Re}{\operatorname{Re}}
\renewcommand{\Im}{\operatorname{Im}}
\newcommand{\wt}{\widetilde}

\def\eqnn{\begin{eqnarray*}}
\def\eeqnn{\end{eqnarray*}}
\def\eqn{\begin{eqnarray}}
\def\eeqn{\end{eqnarray}}
\newcommand{\nc}{\newcommand}
\nc{\be}{\begin{equation}}
\nc{\ee}{\end{equation}}
\nc{\ba}{\begin{eqnarray}}
\nc{\ea}{\end{eqnarray}}
\nc{\eps}{\epsilon}
\def\prf{\begin{proof}}
\def\endprf{\end{proof}}
\begin{document}

\title{Well--posedness and nonlinear smoothing for the ``good" Boussinesq equation on the half-line.}

\author{{\bf E. ~Compaan, N.~Tzirakis}\\
University of Illinois\\
Urbana-Champaign}

\thanks{Email addresses: compaaan2@illinois.edu (E. Compaan), tzirakis@illinois (N. Tzirakis)}\thanks{The first author was supported by a National Physical Science Consortium fellowship. The second author's work was supported by a grant from the Simons Foundation (\#355523 Nikolaos Tzirakis).}

\date{}

\begin{abstract}
In this paper we study the regularity properties of the ``good" Boussinesq equation on the half line. We obtain local existence, uniqueness and continuous dependence on initial data in low-regularity spaces. Moreover we prove that the nonlinear part of the solution on the half line is smoother than the initial data, obtaining half derivative smoothing of the nonlinear term in some cases. Our paper improves the result in \cite{HM}, being the first result that constructs solutions for the initial and boundary value problem of the ``good" Boussinesq equation below the $L^2$ space.  Our theorems are sharp within the framework of the restricted norm method that we use and match the known results on the full line in \cite{KPV} and \cite{F}.
\end{abstract}

\maketitle
\section{Introduction}

We are concerned with the following initial-boundary value problem on the half line, known as the ``good" Boussinesq equation: 
\begin{equation}\label{eq:B} 
 \begin{cases}
  u_{tt} - u_{xx} + u_{xxxx} + (u^2)_{xx} = 0, \quad x \in \R^+, \; t \in \R^+ \\
  u(0,t) = h_1(t),\qquad    u_{x}(0,t) = h_2(t), \\
  u(x,0) = f(x),  \qquad    u_t(x,0) = g_x(x). 
 \end{cases}
\end{equation}
The data $(f,g,h_1,h_2)$ will be taken in the space $H^s_x(\R^+)\times H^{s-1}_x(\R^+) \times H^{\frac{2s + 1}{4}}_t(\R^+) \times H^{\frac{2s-1}{4}}_t(\R^+)$ with the additional compatibility conditions $h_1(0)=f(0)$ when $\frac{1}{2}< s_0\leq\frac32$ and $h_1(0)=f(0),\ \ h_{2}(0)=f^{\prime}(0) $ when $\frac{3}{2}< s_0\leq\frac52$.
These compatibility conditions are necessary since the solutions we are interested in are continuous space-time functions for $s>\frac12$. 

This equation is known as the ``good'' Boussinesq, in contrast to that with the opposite sign in front of the fourth derivative, which was derived by Boussinesq \cite{B} as a water wave model. It also appears as a model of a nonlinear string \cite{Z}. This original Boussinesq equation is linearly unstable because of exponential growth in Fourier modes. The ``good'' Boussinesq \eqref{eq:B} has appeared in studies of shape-memory alloys \cite{FLS}. 
The Boussinesq equation has been extensively studied on $\R$ and $\T$. Bona and Sachs showed well-posedness for data $(f,g) \in H^{s}(\R) \times H^{s-1}(\R)$ for $s > \frac52$ \cite{BS}. Linares established well-posedness for data in $L^2(\R) \times H^{-1}(\R)$ \cite{L} using Strichartz estimates and the theory which Kenig, Ponce and Vega developed for the KdV equation in \cite{KPV1}. Well-posedness in $H^{-\frac14}(\R) \times H^{-\frac54}(\R)$ was shown in \cite{F}, where the restricted norm method of Bourgain ($X^{s,b}$ method, \cite{bou, bou1}) was used. The result in \cite{F} is sharp in the sense that the key bilinear estimate used in the $X^{s,b}$ theory fails for any $s<-\frac{1}{4}$. A simple gauge transformation, \cite{KT}, reduces the ``good" Boussinesq equation into a quadratic nonlinear Schr\"odinger equation, but it is not clear how one can take advantage of this transformation on the half-line. Later in \cite{KT} and \cite{K}, a modification of the restricted norm method of Bourgain was introduced. The well-posedness theory was then improved for both the real line and the torus. In particular, for the real line local well-posedness was established in $H^{-\frac12} \times H^{-\frac32}$. The well-posedness theory at the $H^{-\frac12} \times H^{-\frac32}$ level is known to be sharp, \cite{K}. Our result is sharp, up to an endpoint, in the sense that we also obtain local well-posedness in $H^{-\frac14+}(\R+) \times H^{-\frac54+}(\R+)$, noting that it is not obvious how one can modify the $X^{s,b}$ norm and use an appropriate transformation to simplify the equation in the case of the initial-boundary value problem. 

In this paper we continue the program initiated in \cite{ET} of establishing the regularity properties of nonlinear dispersive partial differential equations (PDE) on a half line using the tools that are available in the case of the real line, where the PDE are fully dispersive. To this end, we extend the data into the whole line and use Laplace transform methods to  set up an equivalent integral equation (on $\R\times \R$) for the solution; see \eqref{eq:decomp} below.  We analyze the integral equation  using the restricted norm method and multilinear $L^2$ convolution estimates. To state the main theorem of this paper we start with a definition.

\begin{definition}
We say that the Boussinesq equation \eqref{eq:B} is locally well-posed in $H^{s}(\R^+)$ if for any $(f,g,h_1,h_2) \in H^s_x(\R^+)\times H^{s-1}_x(\R^+) \times H^{\frac{2s + 1}{4}}_t(\R^+) \times H^{\frac{2s-1}{4}}_t(\R^+)$, with the additional compatibility conditions mentioned above, the equation $\Phi( u) = u$, where $\Phi$ is defined by \eqref{eq:decomp}, has a unique solution in 
\[ X^{s,b}_T \cap C_t^0H^s_x \cap C_x^0H_t^{\frac{2s+1}{4}},\]
for some $b < \frac12$ and some sufficiently small $T$, dependent only on the norms of the initial and boundary data. Furthermore, the solution depends continuously on the initial and boundary data. 
\end{definition}

Our main theorem is below. Note that it extends the result in \cite{HM}, which established well-posedness for $s > \frac12$. In addition we prove that the nonlinear part of the solution is smoother that the initial data. As expected the smoothing disappears at the upper endpoint $s=\frac52$ but not on the lower endpoint $s=-\frac{1}{4}$, where one can still gain a quarter of a derivative. We consider this as an indication (along with the smoothing of order $s+\frac12$) that the ``good" Boussinesq equation should be well-posed in $H^{-\frac12}(\R+) \times H^{-\frac32}(\R+)$, although a modification of our method will be definitely needed to overcome the failure of the bilinear estimates below $H^{-\frac{1}{4}}$. The reader can consult \cite{ETbook} for  many examples of dispersive PDE that enjoy nonlinear smoothing properties at regularities equal to the regularities of the sharp local well-posedness theory. We finally note that the operator  $W_0^t$ is the linear part of the solution of the equation \eqref{eq:B}, see Section 3 below.
\begin{theorem} \label{main_thrm}
 For any $s \in \left(-\frac14, \frac52 \right)$, $s \neq \frac12, \frac32$, the equation \eqref{eq:B} is locally well-posed in $H^s(\R^+)$. Moreover, we have the following smoothing estimate. For $a < \min\{\frac{1}{2} , s + \frac12, \frac52 - s\}$, 
 \begin{equation*}
  u - W_0^t(f,g,h_1,h_2) \in C_0^tH^{s + a}_x.
 \end{equation*} 
 In addition, the solutions are independent of the extensions of the initial data.
\end{theorem}

To prove the above theorems we rely on a Duhamel formulation of the nonlinear system adapted to the boundary conditions. This expresses the nonlinear solution as the superposition of the linear evolutions which incorporate the boundary and the initial data with the nonlinearity. Thus, we first solve two linear problems  by a combination of Fourier and Laplace transforms, \cite{ET}, \cite{bonaetal}, after extending the initial data to the whole line. The idea is then to use the restricted norm method in the Duhamel formula. The uniqueness of the solutions thus constructed is not immediate since we do not know that the fixed points of the Duhamel operators have restrictions on the half line which are independent of the extension of the data. For the case of more regular data the uniqueness property of the solution is proved in \cite{HM}. For less regular data we take advantage of the smoothing estimate we establish in Theorem \ref{main_thrm} to obtain uniqueness all the way down to the local theory threshold $H^{-\frac14+}(\R+) \times H^{-\frac54+}(\R+)$. We remark that this iteration is successful because the full nonlinear estimate we provide remains valid for any $s>-\frac{1}{4}$, matching thus the regularity of the local theory.

As we have already mentioned our result improves the result in \cite{HM}. The initial and boundary value problem (IVBP) for the ``good" Boussinesq equation on the half line has also been considered in \cite{X} and \cite{X2}. In the first paper the author obtained local well-posedness for any $s>\frac{1}{2}$ (having a different set of boundary data than \cite{HM}), while in the second paper the same author obtained local well-posedness for $L^2$ solutions. As far as we know our paper is the first result where well-posed solutions are constructed below the $L^2$ space for the ``good" Boussinesq equation. At this level of regularity, Strichartz type estimates available on the full line are not useful in the construction of solutions obtained through fixed point theorems.

We now discuss briefly the organization of the paper.  In Section \ref{notation}, we introduce some notation and the function spaces that we use to obtain the well-posedness of the IBVP. In Section \ref{results} we define the notion of the solution. More precisely we set up the integral representation (Duhamel's formula) of the nonlinear solution map that we later prove is a contraction in an appropriate metric space. We obtain the solution as a superposition of a linear and a nonlinear evolution. The solution of the linear IBVP can be found by a direct application of the Fourier and the Laplace transform methods. Section \ref{estimates} states the linear and nonlinear a priori estimates that we use to iterate the solution using the restricted norm method appropriately modified for our needs. In Section \ref{main_thrm_prf} we put all the estimates together and show why the solution map is a contraction thus proving the first part of Theorem \ref{main_thrm}. Uniqueness is proved on Section \ref{uniqueness}. Section \ref{Proofs} is the main part of the paper where all the estimates, linear and nonlinear, are established. Finally in Section \ref{appendix} we provide an appendix which justifies the application of the Laplace transform on the half line and the representation formula for the solution of the linear problem with zero initial data.

\section{Notation \& Function Spaces}\label{notation}
We define the one-dimensional Fourier transform by
\[ \widehat{f}(\xi) = \mathcal{F}_x f(\xi) = \int_{\R} e^{-i x \xi} f(x) \d x. \]
We set $\lb \xi \rb = \sqrt{ 1 + |\xi|^2}$. 
The characteristic function on $[0,\infty)$ is denoted by $\chi$.  Sobolev spaces $H^s(\R^+)$ on the half-line for $s > -\frac12$ are defined by 
\begin{align*}
H^s(\R^+) &= \Bigl\{ g \in \mathcal{D}(\R^+) \; : \; \text{there exists } \tilde{g} \in H^s(\R) \text{ with } \tilde{g} \chi = g \Bigr\},  \\
\|g\|_{H^s(\R^+)}&= \inf \Bigl\{ \|\tilde{g} \|_{H^s(\R)} \; : \; \tilde{g}\chi = g \Bigr\}. 
\end{align*}
The restriction $s > -\frac12$ is needed because multiplication with characteristic functions is not defined for $H^s$ distributions when $s \leq -\frac12$. We will also use the $X^{s,b}$ spaces (\cite{bou, bou1}) corresponding to the Boussinesq flow. These are defined for functions on the full space $\R_x \times \R_t$ by the norm 
\[ \|u\|_{X^{s,b}} = \Bigl\| \lb \xi \rb^s \bigl\lb |\tau| - \sqrt{\xi^2 + \xi^4} \bigr\rb^b \widehat{u}(\xi,\tau) \Bigr\|_{L^2_\xi L^2_\tau}. \] 
It is helpful to note (\cite{F}) that there exists $c$ such that 
\[  \frac1c \leq \frac{\bigl\lb a  - \sqrt{b + b^2} \bigr\rb}{\lb a- b \rb} \leq c \quad \text{ for all } \quad a,b \geq 0,\]
so the above $X^{s,b}$ norm is equivalent to $\bigl\| \lb \xi \rb^s \left\lb |\tau| -\xi^2 \right\rb^b \widehat{u}(\xi,\tau) \bigr\|_{L^2_\xi L^2_\tau}$. 

The solution to the linear problem $w_{tt} - w_{xx} + w_{xxxx} = 0$ on $\R$ with initial data $w(x,0) = f(x)$ and $w_t(x,0) = g_x(x)$ will be denoted by
\begin{equation*}
 W_R^t\bigl(f(x),g(x)\bigr) = W_{R,1}^t f(x) + W_{R,2}^t g_x(x), 
\end{equation*}
where $W_{R,1}^t$ and $W_{R,2}^t$ are the Fourier multiplier operators with multipliers $\operatorname{Re} e^{it \sqrt{\xi^2 + \xi^4}}$ and $ \operatorname{Im}e^{it \sqrt{\xi^2 + \xi^4}} (\xi^2 + \xi^4)^{-1/2}$ respectively. 

Let $\rho \in C^\infty$ be a cut-off function such that $\rho = 1$ on $[0, \infty)$ and $\operatorname{supp} \rho \subset [-1, \infty)$. Let $\eta \in C^\infty$ be a bump function such that $\eta = 1$ on $[-1,1]$ and $\operatorname{supp} \eta \subset [-2,2]$. The notation $D_0$ represents evaluation at $x=0$, i.e.
\[ D_0\bigl[ u(x,t) \bigr] = u(0,t). \] 

Finally, the notation $a \lesssim b$ indicates that $a \leq Cb$ for some absolute constant $C$. The expression $a \gtrsim b$ is defined similarly, and $a \approx b$ means that $a \lesssim b$ and $a \gtrsim b$. The notation $a+$ indicates $a + \epsilon$, where $\epsilon$ can be arbitrarily small. We define $a-$ similarly. 

\section{Notion of Solution and Statement of Results} \label{results}

To obtain solutions of \eqref{eq:B}, we begin by constructing the solution of the linear initial-boundary-value problem:
\begin{align} \label{eq:lB}
 \begin{cases}
  v_{tt} - v_{xx} + v_{xxxx}= 0 \\
  v(0,t) = h_1(t), \qquad v_{x}(0,t) = h_2(t), \\
  v(x,0)= f(x), \qquad v_t(x,0) = g(x),
 \end{cases}
\end{align}
with the compatibility condition $h_1(0) = f(0)$ for $\frac12 < s \leq \frac32$, and the additional condition $f'(0) = h_2(0)$ for $\frac32 < s \leq \frac52$.  
Denote this solution by $W^t_0(f,g,h_1, h_2)$. For extensions $f^e$ and $g^e$ to the full line $\R$ of the functions $f$ and $g$, we may write
\[ W_0^t(f, g, h_1, h_2) = W_0^t(0,0, h_1 - p_1, h_2 - p_2) + W_R^t(f^e,g^e), \]
where $p_1(t) =  D_0 \bigl[ W_R^t(f^e,g^e)\bigr]$ and $p_2(t) =  D_0\bigl[ W_R^t(f_e,g_e) \bigr]_x$. We thus decompose the solution operator as a sum of a modified boundary operator, which incorporates zero initial data, and the free propagator defined on the whole real line. For $x>0$, this solution formula expresses the unique solution of \eqref{eq:lB}. Note that $W_0^t(0,0,h_1,h_2)$ is the solution to the following problem: 

\begin{align} \label{eq:zB}
 \begin{cases}
  v_{tt} - v_{xx} + v_{xxxx} = 0 \\
  v(0,t) = h_1(t), \qquad v_{x}(0,t) = h_2(t), \\
  v(x,0)= 0, \quad \qquad v_t(x,0) = 0.
 \end{cases}
\end{align}

We will use the following explicit representation of $W_0^t(0,0,h_1,h_2)$ extensively. It is proved in the Appendix using a Laplace transform argument. Similar expressions have been derived in \cite{X, X2} using the Laplace transform and in \cite{HM} using Fokas' unified transform method.

\begin{lemma}\label{explicit_sol}
Suppose $h_1$ and $h_2$ are Schwarz functions. The solution to \eqref{eq:zB} on $\R^+ \times \R^+$ can be written in the form $v(x,t) = \frac1{2\pi} (-A - B + C + D)$, where 
\begin{align} \label{eq:lin_ibvp}
\begin{split}
A = &\int_{-\infty}^{ \infty} \frac{e^{i t \omega \sqrt{\omega^2 + 1} -x \sqrt{\omega^2 + 1}}}{\sqrt{1 + \omega^2}}i \omega\biggl(i \omega + \sqrt{1 + \omega^2}\biggr)\; \widehat{h_1}\left( \omega \sqrt{\omega^2 +1}\right) \rho\Bigl(x \sqrt{\omega^2 + 1} \Bigr) \d \omega \\
B = & \int_{-\infty}^{ \infty} \frac{e^{i t \omega \sqrt{\omega^2 + 1} - x\sqrt{\omega^2 + 1}}}{\sqrt{1 + \omega^2}}     \biggl(i \omega + \sqrt{1 + \omega^2}\biggr)\; \widehat{h_2}\left( \omega \sqrt{\omega^2 +1}\right) \rho\Bigl(x \sqrt{\omega^2 + 1} \Bigr) \d \omega \\
C = &\int_{-\infty}^{ \infty} e^{i t \omega \sqrt{\omega^2 + 1} - ix\omega}                              \biggl(i \omega + \sqrt{1 + \omega^2}\biggr)            \; \widehat{h_1}\left( \omega \sqrt{\omega^2 +1}\right) \d \omega \\
D = &\int_{-\infty}^{ \infty} \frac{e^{i t \omega \sqrt{\omega^2 + 1} - ix\omega}}{\sqrt{1 + \omega^2}}     \biggl(i \omega + \sqrt{1 + \omega^2}\biggr)            \; \widehat{h_2}\left( \omega \sqrt{\omega^2 +1}\right) \d \omega. 
\end{split}
\end{align}
Here by an abuse of notation, $\widehat{h_i}$ denotes the Fourier transform of $\chi h_i$. 
\end{lemma}
This explicit form will be used to establish bounds on $W_0^t(0,0,h_1,h_2)$ in the subsequent sections. Notice that the integrals $A$, $B$, $C$, and $D$ are defined on the entire space $\R_x \times \R_t$ thanks to the inclusion of the cut-off function $\rho$. 

It is now clear that the solution to the full initial-boundary-value problem \eqref{eq:B} satisfies, for $t\leq T$, the equation $\Phi(u) = u$, where the operator $\Phi$ is given by
\begin{equation} \label{eq:decomp}
 \begin{split}
 \Phi \bigl(u(x,t)\bigr) = \eta(t/T) W_{R}^t \bigl(f^e(x),g^e(x) \bigr) + \eta(t/T) \int_0^t W_{R,2}^{t-t'} G(u) \d t' \\
  + \eta(t/T) W^t_0(0,0,h_1 - p_1 - q_1, h_2 - p_2 - q_2),
 \end{split}
\end{equation}
with
\begin{equation}\label{eq:nonlinearity} 
G(u) = \eta(t/T) (u^2)_{xx}, 
\end{equation}
\begin{equation} \label{eq:forcing_terms}
\begin{split}
 p_1(t) = \eta(t/T) D_0\Bigl[W_{R}^t \bigl(f^e(x),g^e_x(x) \bigr)\Bigr],     \quad q_1(t) = \eta(t/T) D_0 \Bigl[ \int_0^t W_{R,2}^{t-t'} G(u) \d t'\Bigr], \\
 p_2(t) = \eta(t/T) D_0\Bigl[W_{R}^t \bigl(f^e(x),g^e_x(x) \bigr)\Bigr]_{x}, \quad q_2(t) = \eta(t/T) D_0 \Bigl[ \int_0^t W_{R,2}^{t-t'} G(u) \d t' \Bigr]_{x}.
 \end{split}
\end{equation}

In the following, we will use a fixed point argument to obtain a unique solution to $\Phi( u )= u$ in a suitable function space on $\R \times \R$ for sufficiently small $T$. The restriction of $u$ to $\R^+ \times \R$ is a distributional solution of \eqref{eq:B}. Furthermore, smooth solutions of $\Phi(u) = u$ are classical solutions of  \eqref{eq:B}. 

The contraction argument is carried out in $X^{s,b}$ spaces. To bound the solution to the linear Boussinesq on $\R$ and the Duhamel term, we will use the following estimates from \cite{F}. 
For any $s$ and $b$, we have 
\begin{equation} \label{eq:lin_est}
 \| \eta(t) W_{R}^t \bigl( f, g \bigr)\|_{X^{s,b}} \lesssim \| f\|_{H^s} + \|g\|_{H^{s-1}}. 
\end{equation}
Furthermore, for any $-\frac12 < b' \leq 0 \leq b \leq b' + 1$ and $0 < T< 1$, the estimate 
\begin{equation} \label{eq:duhamel_est}
\Bigl\| \eta(t/T) \int_0^t W_{R,2}^{t-t'} G(u) \d t' \Bigr\|_{X^{s,b}} \lesssim T^{1-(b-b')} \Bigl\|\mathcal{M}(G(u)) \Bigr\|_{X^{s, b'}} 
\end{equation}
holds, where $\mathcal{M}$ is the Fourier multiplier operator defined by $\widehat{\mathcal{M}(f)} = (\xi^2 + \xi^4)^{-1/2} \widehat{f}$. Also 
\begin{equation}\label{eq:t_power}
\| \eta(t/T) F\|_{X^{s,b_1}} \lesssim T^{b_2 - b_1} \| F\|_{X^{s,b_2}},
\end{equation}
for any $-\frac12 < b_1 < b_2 < \frac12$ \cite{ETbook}. Finally, we require the following lemma regarding extensions of $H^s(\R^+)$ functions. It will be used to bound the explicit linear solution given in Lemma \ref{explicit_sol}, which is given in terms of the Fourier transforms of $\chi h_i$.  
\begin{lemma}{\cite{ET}} \label{char_func_lemma} Assume $h \in H^s(\R^+)$. 
\begin{enumerate}
 \item If $-\frac12 < s < \frac12$, then $\| \chi h \|_{H^s(\R)} \lesssim \| h \|_{H^s(\R^+)}$. 
 \item If $ \frac12 < s < \frac32$ and $h(0) = 0$, then $\| \chi h \|_{H^s(\R)} \lesssim \| h \|_{H^s(\R^+)}$. 
\end{enumerate}
\end{lemma}

\section{A Priori Estimates}\label{estimates}

To close the contraction argument, we need a number of estimates on the terms in \eqref{eq:decomp}.

\subsection{Linear Estimates}

First, we give a Kato smoothing inequality, which is proved in Section \ref{Kato_prf}. Similar results are stated in \cite{X}. This estimate is necessary to ensure that $\Phi(u)$ lies in $L^\infty_xH_t^{\frac{2s+1}{4}}$ and to control the terms $p_i$ defined in \eqref{eq:forcing_terms}. 
\begin{lemma}\label{Kato}
 For any $s$, 
 \begin{align*}
  \| \eta(t) W_R^t(f,g) \|_{{L^\infty_x}H_t^{\frac{2s+1}{4}}} \lesssim \|f\|_{H^s_x} + \|g\|_{H^{s-1}_x} \\
  \| \eta(t) [W_R^t(f,g)]_{x} \|_{{L^\infty_x}H_t^{\frac{2s-1}{4}}} \lesssim \| f\|_{H^s_x} + \|g\|_{H^{s-1}_x}. 
 \end{align*} 
\end{lemma}

For the solution to the linear initial-boundary-value problem we have the following estimates, which are proved in Sections \ref{linear_half_prf} and \ref{cont_prf}. These are used to bounded the $W_0^t$ term in $\Phi(u)$, and to ensure that this term lies in the desired space $C_t^0H_x^s \cap C_x^0 H_t^{\frac{2s+1}{4}}$. 

 \begin{lemma}\label{linear_half}
For any compactly supported smooth function $\eta$ and any $s \geq - \frac12$ with $b < \frac12$, 
 \[ \| \eta(t) W^t_0(0,0,h_1, h_2) \|_{X^{s,b}} \lesssim \| \chi h_1 \|_{H_t^{\frac{2s+1}{4}}(\R)} + \|\chi  h_2 \|_{H_t^{\frac{2s-1}{4}}(\R)}.\] 
\end{lemma}

\begin{lemma}\label{cont}
For any $s \geq -1 $ and intial data $(h_1, h_2)$ such that $(\chi h_1,\chi h_2) \in H^{\frac{2s+1}{4}}(\R) \times H^{\frac{2s-1}{4}}(\R)$, we have 
\begin{align*}
 W_0^t(0,0,h_1,h_2) &\in C_t^0H_x^s(\R \times \R) \\
 \eta(t) W_0^t(0,0,h_1,h_2) &\in C_x^0 H_t^{\frac{2s+1}{4}}( \R \times \R). 
 \end{align*}
\end{lemma}

\subsection{Nonlinear Estimates}
\begin{lemma}\label{bilin}
Let $\mathcal{M}$ be the Fourier multiplier operator with multiplier $(\xi^2 + \xi^4)^{-1/2}$. For $s > - \frac14$ with $a < \min\{\frac12, s + \frac12\}$ and $\frac12 - b > 0 $ sufficiently small, we have 
 \[ \|\mathcal{M}(uv)_{xx} \|_{X^{s+a, -b}} \lesssim \| u\|_{X^{s,b}}\| v\|_{X^{s,b}}. \]
\end{lemma}

This lemma is proved in Section \ref{bilin_prf}. The next requirement is to control the Duhamel part of the correction term. This is accomplished with the following estimate, which is proved in Section \ref{Duhamel_Kato_prf}. 
\begin{lemma}\label{Duhamel_Kato} 
For $\frac12 - b > 0$ sufficiently small, we have 
\begin{align*}
\left\| \eta(t) \int_0^t W_{R,2}^{t-t'} G \d t' \right\|_{L^\infty_xH^{\frac{2s+ 1}{4}}_t} + \left\| \eta(t) \Bigl[ \int_0^t W_{R,2}^{t-t'} G \d t'\Bigr]_{x} \right\|_{L^\infty_xH^{\frac{2s- 1}{4}}_t}\hspace{1.5in} \\
\lesssim 
 \begin{cases} 
 \|\mathcal{M}(G)\|_{X^{s,-b}} + \Bigl\| \lb \tau \rb^{\frac{2s-1}4} \int \chi_{Q}(\xi,\tau) \lb \xi \rb^{-1}| \widehat{\mathcal{M}(G)}(\xi,\tau)|\d \xi \Bigr\|_{L^2_\tau}&\text{ if } -\frac12 \leq s \leq \frac12 \\
 \|\mathcal{M}(G)\|_{X^{s,-b}} + \Bigl\| \int \chi_R(\xi, \tau)  \lb |\tau| - \xi^2 \rb^\frac{2s-3}{4} | \widehat{\mathcal{M}(G)}(\xi,\tau)| \d \xi \Bigr\|_{L^2_\tau} &\text{ if } s > \frac12,
 \end{cases} 
\end{align*}
where $Q = \{ \| \tau | \ll \xi^2\} \cap \{ |\xi| \gtrsim 1\}$ and $R = \{ |\tau| \gg \xi^2 \} \cup \{ |\xi| \lesssim 1\}$. 
\end{lemma}

It remains to bound the left hand side of the inequality in Lemma \ref{Duhamel_Kato}. We use Lemma \ref{bilin} to control the $X^{s,b}$ norms; the other terms are bounded using the following lemmata, which are proved in Sections \ref{correction_prf} and \ref{bilinear_prf} respectively.  
\begin{lemma}\label{bilinear0} 
Let $Q$ be the set $\{ |\tau| \ll \xi^2\} \cap \{ |\xi| \gtrsim 1\}$. For $ -\frac14 < s + a \leq \frac12$ and $0 \leq a < s + \frac12$, we have 
\[ \Bigl\| \lb \tau \rb^{\frac{2(s+a)-1}{4}} \int \chi_Q(\xi,\tau) \lb \xi \rb^{-1} \frac{\xi^2}{\sqrt{\xi^2 + \xi^4}} |\widehat{uv}(\xi,\tau)| \d \xi \Bigr\|_{L^2_\tau} \lesssim \|u\|_{X^{s,b}}\|v\|_{X^{s,b}}.  \]
\end{lemma}

\begin{lemma}\label{bilinear} 
Let $R$ be the set $\{ |\tau| \gg \xi^2\} \cup \{ |\xi| \lesssim 1\}$. For $ \frac12 < s + a \leq \frac52$ and $a < \min\{1, s + \frac12\}$, we have 
\[ \Bigl\| \int \chi_R(\xi,\tau) \lb |\tau| -  \xi^2 \rb^\frac{2(s+a)-3}{4} \frac{\xi^2}{\sqrt{\xi^2 + \xi^4}} |\widehat{uv}(\xi,\tau)| \d \xi \Bigr\|_{L^2_\tau} \lesssim \|u\|_{X^{s,b}}\|v\|_{X^{s,b}}.  \]
\end{lemma}

\section{Local Theory: Proof of Theorem \ref{main_thrm}}\label{main_thrm_prf}

We will first show that the map $\Phi$ defined in \eqref{eq:decomp} has a unique fixed point in $X^{s,b}$. 
Let $f^e \in H^s(\R)$ and $g^e \in H^{s-1}(\R)$ be extensions of $f$ and $g$ such that $\|f^e\|_{H^s(\R)} \lesssim \|f\|_{H^s(\R^+)}$ and $\|g^e\|_{H^{s-1}(\R)} \lesssim \|g\|_{H^{s-1}(\R^+)}$. Recall that  
\begin{equation}\label{Duhamel}
 \begin{split}
 \Phi \bigl(u(x,t)\bigr) = \eta(t/T) W_{R}^t \bigl(f^e(x),g^e(x) \bigr) + \eta(t/T) \int_0^t W_{R,2}^{t-t'} G(u) \d t' \\
  + \eta(t/T) W^t_0\bigl(0,0,h_1 - p_1 - q_1, h_2 - p_2 - q_2\bigr),
 \end{split}
\end{equation}
where $G(u)$, $p_i$, and $q_i$ are defined in \eqref{eq:nonlinearity}-\eqref{eq:forcing_terms}. 
To bound the first summand in $\Phi$, apply \eqref{eq:lin_est} to obtain  
\[ \left\| \eta(t/T) W_{R}^t \bigl(f^e,g^e \bigr) \right\|_{X^{s,b}} \lesssim \|f^e\|_{H^s(\R)} + \|g^e\|_{H^{s-1}(\R)} \lesssim \|f\|_{H^s(\R^+)} + \|g\|_{H^{s-1}(\R^+)}. \]
For the Duhamel term, we apply \eqref{eq:duhamel_est} and Lemma \ref{bilin} to obtain 
\begin{align*}
 \left\| \eta(t/T) \int_0^t W_{R,2}^{t-t'} G(u) \d t' \right\|_{X^{s,b}} \lesssim T^{1-2b} \| \mathcal{M}\bigl(u^2\bigr)_{xx} \|_{ X^{s,-b}} \lesssim  T^{1-2b} \| u\|^2_{X^{s,b}}. 
\end{align*}
Finally, for the $W_0^t$ term, we apply Lemma \ref{linear_half} and Lemma \ref{char_func_lemma} to obtain
\begin{align*}
 &\left\| \eta(t/T) W^t_0\bigl(0,0,h_1 - p_1 - q_1, h_2 - p_2 - q_2\bigr) \right\|_{X^{s,b}} \\
 &\lesssim \| \chi(h_1 -p_1 -q_1) \|_{H_t^{\frac{2s+1}{4}}(\R)} + \| \chi(h_2 -p_2 -q_2) \|_{H_t^{\frac{2s-1}{4}}(\R)} \\
 &\lesssim \|h_1 -p_1 \|_{H_t^{\frac{2s+1}{4}}(\R^+)} + \| q_1 \|_{H_t^{\frac{2s+1}{4}}(\R^+)} + \| h_2 -p_2\|_{H_t^{\frac{2s-1}{4}}(\R^+)} + \|q_2\|_{H_t^{\frac{2s-1}{4}}(\R^+)}.
\end{align*}
By Kato smoothing, Lemma \ref{Kato}, we have 
\begin{align*}
 \|p_1 \|_{H_t^{\frac{2s+1}{4}}(\R)} + \| p_2\|_{H_t^{\frac{2s-1}{4}}(\R)} \lesssim \| f^e\|_{H^s(\R)} + \|g^e\|_{H^{s-1}(\R)} \lesssim \|f\|_{H^s(\R^+)} + \|g\|_{H^{s-1}(\R^+)}.
\end{align*}
To bound the $q_i$ norms, we apply Lemma \ref{Duhamel_Kato}, \eqref{eq:t_power}, and Lemma \ref{bilin}, Lemma \ref{bilinear0}, and Lemma \ref{bilinear} to obtain the bounds
\[ \| q_1\|_{H_t^{\frac{2s+1}{4}}(\R)} +  \| q_2\|_{H_t^{\frac{2s+1}{4}}(\R)} \lesssim T^{\frac12 - b - } \| u\|_{X^{s,b}}^2. \]
Combining these estimates, we find that 
\[ \left\| \Phi(u) \right\|_{X^{s,b}} \lesssim  \|f\|_{H^s(\R^+)} + \|g\|_{H^{s-1}(\R^+)} + \|h_1 \|_{H_t^{\frac{2s+1}{4}}(\R^+)} +  \| h_2 \|_{H_t^{\frac{2s-1}{4}}(\R^+)} + T^{\frac12 - b - } \| u\|_{X^{s,b}}^2.\] 

This, together with similar estimates for the difference $\Phi(u) - \Phi(v)$, yields the existence of a fixed point of $\Phi$ for $T$ sufficiently small: 
\[T = T\bigl(\|f\|_{H^s(\R^+)}, \|g\|_{H^{s-1}(\R^+)}, \|h_1 \|_{H_t^{\frac{2s+1}{4}}(\R^+)},  \| h_2 \|_{H_t^{\frac{2s-1}{4}}(\R^+)}\bigr).\] 

Next, we establish continuity in $H^s$. For the $W_0^t$ term, this follows from Lemma \ref{cont}. The first term of $\Phi$, the linear flow on $\R$, can be seen to be continuous from its Fourier multiplier formula. Continuity of the Duhamel term follows from the embedding $X^{s,b} \subset C^0_tH^{s}_x$ for $b > \frac12$ along with \eqref{eq:duhamel_est} and Lemma \ref{bilin}. The fact that the solution lies in $C^0_xH^{\frac{2s+1}{4}}_t$ follows from Lemma \ref{cont} for the $W_0^t$ term, from Kato smoothing (Lemma \ref{Kato}) for the linear flow on $\R$, and from Lemmata \ref{Duhamel_Kato}-\ref{bilinear} for the Duhamel term. 

\section{Uniqueness of Solutions}\label{uniqueness}

In this section, we show that solutions to \eqref{eq:B} derived in the previous section are unique. For $s > \frac12$, uniqueness of $C^0([0,T], H^s(\R^+))$ solutions to \eqref{eq:B} holds by \cite{HM}. The solutions obtained in the previous section also lie in this space after restriction to $x \in \R^+$. Thus we have uniqueness for $s > \frac12$.

Using the smoothing estimates in Theorem \ref{main_thrm}, we can now obtain uniqueness of local solutions for the full range of Sobolev exponents in the local theory. First consider initial data $(f,g,h_1, h_2) \in H^s_x(\R^+)\times H^{s-1}_x(\R^+) \times H^{\frac{2s + 1}{4}}_t(\R^+) \times H^{\frac{2s-1}{4}}_t(\R^+)$ for some $s \in \bigl(0, \frac12\bigr)$. Suppose $f^e$ and $\wt{f}^e$ are two $H^{s}(\R)$ extensions of $f$, and $g^e$ and $\wt{g}^e$ are two $H^{s-1}(\R)$ extensions of $g$. Let $u$ and $\wt{u}$ be the corresponding solutions of the fixed-point equation for $\Phi$. Take a sequence $f_k \in H^{\frac12+}(\R^+)$ converging to $f$ in $H^s(\R^+)$. Let $f^e_k$ and $\wt{f}^e_k$ be $H^{\frac12+}(\R)$ extensions of $f_k$ which converge to $f^e$ and $\wt{f}^e$ respectively in $H^r(\R)$ for $r < \frac12-$. This is possible by Lemma \ref{approx_lemma} below. Similarly, obtain a sequence $g_k$ in $H^{-\frac12+}(\R^+)$ converging to $g$ in $H^{s-1}(\R^+)$, and extensions $g^e_k$ and $\wt{g}^e_k$ of $g_k$ converging to $g^e$ and $\wt{g}^e$ respectively. 

Using a contraction argument on the set
\begin{align*}
 &\Bigl\{ u \; : \; \|u\|_{X^{\frac12+,b}} \leq C \bigl(\|f_k\|_{H^{1/2+}_x(\R^+)}+ \|g_k\|_{H^{-1/2+}_x(\R^+)} + \|h_1\|_{H^{1/2+}_t(\R^+)} + \| h_2\|_{H^{0+}_t(\R^+)} \bigr) \Bigr\} \\
  & \cap \; \Bigl\{ u \; : \; \|u\|_{X^{s,b}} \leq C \bigl(\|f\|_{H^{s}_x(\R^+)}+ \|g\|_{H^{s-1}_x(\R^+)} + \|h_1\|_{H^{\frac{2s+1}{4}}_t(\R^+)} + \| h_2\|_{H^{\frac{2s-1}{4}}_t(\R^+)} \bigr) \Bigr\},
\end{align*}
we construct $H^{\frac12+}$ solutions of $u_k$ and $\wt{u}_k$ of the Boussinesq \eqref{eq:B} using the extensions $(f^e_k, g^e_k)$ and $(\wt{f}^e_k, \wt{g}^e_k)$ respectively. The smoothing estimates give us a time of existence proportional to the data in the lower norm: $T= T\bigl(\|f\|_{H^s(\R^+)}, \|g\|_{H^{s-1}(\R^+)}, \|h_1 \|_{H_t^{\frac{2s+1}{4}}(\R^+)},  \| h_2 \|_{H_t^{\frac{2s-1}{4}}(\R^+)}\bigr)$. By uniqueness of $H^{\frac12+}$ solutions, $u_k$ and $\wt{u}_k$ are equal on $\R^+ \times \R^+$. By the fixed-point argument, $u_k$ and $\wt{u}_k$ converge in $H^{s-}$ to $u$ and $\wt{u}$ respectively. Thus $u = \wt{u}$ on $\R^+ \times \R^+$. Iterating this argument, we obtain uniqueness for $s > -\frac14$. 

\begin{lemma}{\cite{ET1}}\label{approx_lemma}
 Fix $-\frac12 < s < \frac12$ and $k > s$. Let $p \in H^s(\R^+)$ and $q \in H^k(\R^+)$. Let $p^e$ be an $H^s$ extension of $p$ to $\R$. Then there is an $H^k$ extension $q^e$ of $q$ to $\R$ such that 
 \[ \| p^e - q^e \|_{H^r(\R)} \lesssim \| p-q\|_{H^s(\R^+)} \quad \text{ for } r < s. \] 
\end{lemma}

\section{Proofs of Estimates} \label{Proofs}

\subsection{Proof of Lemma \ref{Kato}: Kato Smoothing} \label{Kato_prf}

We wish to show that
\begin{align*}
  \| \eta(t) W_R^t(f,g) \|_{{L^\infty_x}H_t^{\frac{2s+1}{4}}} \lesssim \|f\|_{H^s_x} + \|g\|_{H^{s-1}_x} \\
  \| \eta(t) [W_R^t(f,g)]_{x} \|_{{L^\infty_x}H_t^{\frac{2s-1}{4}}} \lesssim \| f\|_{H^s_x} + \|g\|_{H^{s-1}_x}. 
 \end{align*} 

 It suffices to consider evaluation at $x = 0$ since Sobolev norms are invariant under translations. Using the Fourier multiplier form of the linear flow, write
 \begin{align*}
  2&\mathcal{F}_t\bigl(\eta W_R^t(f,g)\bigr)(0,\tau) = \int \widehat{\eta}(\tau - \sqrt{\omega^2 + \omega^4}) \;\widehat{f}(\omega) \d \omega + \int \widehat{\eta}(\tau + \sqrt{\omega^2 + \omega^4})\; \widehat{f}(\omega) \d \omega \\
  &+ \int \widehat{\eta}(\tau - \sqrt{\omega^2 + \omega^4}) \frac{\omega}{\sqrt{\omega^2 + \omega^4}}  \;\widehat{g}(\omega) \d \omega - \int \widehat{\eta}(\tau + \sqrt{\omega^2 + \omega^4})  \frac{\omega}{\sqrt{\omega^2 + \omega^4}}\; \widehat{g}(\omega) \d \omega.
 \end{align*}
On the region where $\omega \leq 1$, these terms can easily be bounded in $H_t^{\frac{2s+1}{4}}$ since $\eta$ is a Schwarz function. When $|\omega| > 1$, change variables by setting $\lambda = \omega \sqrt{\omega^2 +1}$. The first two integrals in the above sum are then of the form
\[ \int_{|\lambda | > 2} \widehat{\eta}(\tau \pm |\lambda|) \;\widehat{f}(\omega(\lambda)) \frac{\sqrt{1 + \omega(\lambda)^2}}{2 \omega(\lambda)^2 + 1} \d \lambda,\]
and we wish to bound 
\begin{align*}
 \left\| \int_{|\lambda | > 2} \lb \tau \rb^{\frac{2s + 1}{4}} \widehat{\eta}(\tau \pm |\lambda|) \;\widehat{f}(\omega(\lambda)) \frac{\sqrt{1 + \omega(\lambda)^2}}{2 \omega(\lambda)^2 + 1} \d \lambda \right\|_{L^2_\tau}. 
\end{align*}
Note that the inequality $\lb a + b \rb \lesssim \lb a \rb \lb b \rb$ implies that for any $\alpha$, we have $\lb a + b \rb^ \alpha \lesssim \lb a \rb^{|\alpha|} \lb b \rb^\alpha$. Using this, the quantity above is bounded by
\begin{align*}
 \left\| \int_{|\lambda | > 2} \lb \tau \pm |\lambda| \rb^{\frac{|2s + 1|}{4}} \widehat{\eta}(\tau \pm |\lambda|) \lb\lambda \rb^{\frac{2s + 1}{4}} \;\widehat{f}(\omega(\lambda)) \frac{\sqrt{1 + \omega(\lambda)^2}}{2 \omega(\lambda)^2 + 1} \d \lambda \right\|_{L^2_\tau}. 
\end{align*}
Since $\widehat{\eta}$ is a Schwarz function, we may use Young's inequality and then change variables back to $\omega$ to bound this quantity by  
\begin{align*}
  \left\| \lb\lambda \rb^{\frac{2s + 1}{4}} \;\widehat{f}(\omega(\lambda)) \frac{\sqrt{1 + \omega(\lambda)^2}}{2 \omega(\lambda)^2 + 1} \right\|_{L^2_{|\lambda| > 2}} \lesssim \; \|f\|_{H^s(\R)},
\end{align*}
as desired. The remaining integrals, those involving $g$, can be treated in exactly the same way and bounded by $\|g\|_{H^{s-1}(\R)}$. We obtain the bound on $\| \eta(t) [W_R^t(f,g)]_{x} \|_{{L^\infty_x}H_t^{\frac{2s-1}{4}}}$ by the same argument. 

\subsection{Proof of Lemma \ref{linear_half}: Bounds on Linear Solution}\label{linear_half_prf} 
Recall that we wish to establish
\[ \| \eta(t) W^t_0(0,0,h_1, h_2) \|_{X^{s,b}} \lesssim \| \chi h_1 \|_{H_t^{\frac{2s+1}{4}}(\R)} + \|\chi  h_2 \|_{H_t^{\frac{2s-1}{4}}(\R)}, \] 
where $2 \pi\, W^t_0(0,0,h_1, h_2) = -A - B + C +D$, and the terms $A$, $B$, $C$, and $D$ are given in \eqref{eq:lin_ibvp}. Notice that 
\begin{align*} C &= L^t\phi_C, \quad &\text{where} \quad \widehat{\phi_C}(\omega) =  \Bigl( i \omega + \sqrt{\omega^2 + 1} \Bigr) \; \widehat{h_1}\Bigl(\omega\sqrt{\omega^2 +1}\Bigr), \\
D &= L^t\phi_D, \quad &\text{where} \quad \widehat{\phi_D}(\omega) = \frac{i \omega + \sqrt{\omega^2 + 1}}{\sqrt{\omega^2 + 1}} \; \widehat{h_2}\Bigl(\omega\sqrt{\omega^2 +1}\Bigr),
\end{align*}
and $L^t$ is the spatial Fourier multiplier operator with multiplier $e^{it \omega\sqrt{1+\omega^2}}$. The proof of \eqref{eq:lin_est} implies that 
$\|\eta(t)C \|_{X^{s,b}} \lesssim \| \phi_C \|_{H^{s}_x}$  and $\|\eta(t)D \|_{X^{s,b}} \lesssim \| \phi_D \|_{H^{s}_x}$. Now 
\begin{align*}
\| \phi_C \|_{H^s_x}^2 &= \int_{-\infty}^\infty (2\omega^2 + 1 ) \lb \omega \rb^{2s} \; \left| \widehat{h_1}\Bigl(\omega \sqrt{\omega^2 + 1}\Bigr)\right|^2 \d \omega \\
&= \int_{-\infty}^\infty  \frac{(\omega^2 + 1)^{1/2} \lb \omega\rb^{2s}}{ \bigl\lb \omega \sqrt{\omega^2 + 1} \bigr\rb^{\frac{2s + 1}{2}} } \; \bigl\lb \omega \sqrt{\omega^2 + 1} \bigr\rb^{\frac{2s + 1}{2}} \; \left|\widehat{h_1}\Bigl(\omega \sqrt{\omega^2 + 1}\Bigr)\right|^2 \frac{2\omega^2 + 1}{\sqrt{\omega^2 + 1}} \d \omega  \\
&\lesssim \int_{-\infty}^\infty \bigl\lb \omega \sqrt{\omega^2 + 1} \bigr\rb^{\frac{2s + 1}{2}} \; \left|\widehat{h_1}\Bigl(\omega \sqrt{\omega^2 + 1}\Bigr)\right|^2 \frac{2\omega^2 + 1}{\sqrt{\omega^2 + 1}} \d \omega  \\
&= \int_{-\infty}^\infty \lb z \rb^{\frac{2s + 1}{2}} \; \left|\widehat{h_1}(z)\right|^2 \d z = \| \chi h_1\|_{H^{\frac{2s+1}{4}}(\R)},
\end{align*}
where we used the change of variable $z = \omega \sqrt{\omega^2 + 1}$. 
Similarly, 
\begin{align*}
\| \phi_D\|_{H^s_x}^2 &= \int_{-\infty}^\infty \frac{(2\omega^2 + 1)\lb \omega \rb^{2s}}{\omega^2 + 1} \left|\widehat{h_2}\Bigl(\omega \sqrt{\omega^2 + 1}\Bigr)\right|^2 \d \omega \\
&= \int_{-\infty}^\infty \frac{\lb \omega \rb^{2s}}{(\omega^2 + 1)^{1/2} \bigl\lb \omega \sqrt{\omega^2 + 1} \bigr\rb^{\frac{2s - 1}{2}}} \bigl\lb \omega \sqrt{\omega^2 + 1} \bigr\rb^{\frac{2s - 1}{2}}\; \left|\widehat{h_2}\Bigl(\omega \sqrt{\omega^2 + 1}\Bigr)\right|^2 \frac{2\omega^2 + 1}{\sqrt{\omega^2 + 1}} \d \omega \\
& \lesssim \int_{-\infty}^\infty \lb z \rb^{\frac{2s -1}{2}} \; \left|\widehat{h_2}(z)\right|^2 \d z = \|\chi h_2\|_{H^{\frac{2s-1}{4}}(\R)}.
\end{align*}
Thus we have the desired bounds on $C$ and $D$. 
Now we move on to $A$ and $B$. Assume first that $s = 0$ and $b = \frac12-$. Let $f(y) = e^{-y} \rho(y)$. Then 
\begin{align*}
A(x,t) &= \int_{-\infty}^\infty \frac{i\omega\bigl(i\omega + \sqrt{1 + \omega^2}\bigr)}{\sqrt{\omega^2 + 1}} e^{it \omega \sqrt{\omega^2 + 1}}\; \widehat{h_1}\Bigl( \omega \sqrt{\omega^2 + 1}\Bigr) f\Bigl(x \sqrt{\omega^2 + 1}\Bigr) \d \omega, \\
\widehat{ \eta A}(\xi, \tau) &= \int_{-\infty}^\infty \widehat{\eta}(\tau - \omega\sqrt{\omega^2 + 1})  \frac{i\omega\bigl(i\omega + \sqrt{1 + \omega^2}\bigr)}{\sqrt{\omega^2 + 1}}\; \widehat{h_1}(\omega\sqrt{\omega^2 +1}) \mathcal{F}_x\Bigl(f(x\sqrt{\omega^2 + 1})\Bigr)(\xi)\d \omega \\
&= \int_{-\infty}^\infty \widehat{\eta}(\tau - \omega\sqrt{\omega^2 + 1})  \frac{i\omega\bigl(i\omega + \sqrt{1 + \omega^2}\bigr)}{\omega^2 + 1} \;\widehat{h_1}(\omega\sqrt{\omega^2 +1})\widehat{f}\left(\xi/\sqrt{\omega^2 + 1}\right) \d \omega.
\end{align*}
Since $f$ is a Schwarz function, we have 
\[ |\widehat{f}(\xi/\sqrt{\omega^2 + 1})| \lesssim \frac{1}{1 + \xi^2/(\omega^2 + 1)} = \frac{ \omega^2+1}{1 + \omega^2 + \xi^2}. \] 
Note also that since $\eta$ is a Schwarz is function,
\begin{align*}
| \widehat{\eta}(\tau - \omega\sqrt{\omega^2 + 1})| &\lesssim \lb \tau - \omega \sqrt{\omega^2 + 1} \rb^{-5/2+} \\
& \lesssim \lb \tau - \omega \sqrt{\omega^2 + 1} \rb^{-2} \lb \tau - \xi^2 \rb^{-1/2+} \lb \omega \sqrt{\omega^2 + 1} - \xi^2 \rb^{1/2-}. 
\end{align*}
Therefore, using the bounds for $f$, those for $\eta$, and then moving the $\xi$ norm inside the integral, 
\begin{align*}
\| \eta A\|_{X^{0, \frac12-}} &\lesssim \left\| \lb \tau - \xi^2 \rb^{1/2-} \int _{-\infty}^\infty \left| \widehat{\eta}(\tau - \omega \sqrt{\omega^2 + 1})\right| \frac{|\omega| \sqrt{1 + 2 \omega^2}}{1 + \omega^2 + \xi^2} \left| \widehat{h_1}(\omega\sqrt{\omega^2 + 1})\right| \d \omega \right\|_{L^2_{\xi,\tau}} \\
&\lesssim \left\| \int _{-\infty}^\infty \lb \tau - \omega \sqrt{\omega^2 + 1}\rb^{-2} \frac{|\omega|\sqrt{1 + 2 \omega^2}}{(1 + \omega^2 + \xi^2)^{1/2+}} \left| \widehat{h_1}(\omega\sqrt{\omega^2 + 1})\right| \d \omega \right\|_{L^2_{\xi,\tau}} \\
&\lesssim  \left\| \int _{-\infty}^\infty \lb \tau - \omega \sqrt{\omega^2 + 1}\rb^{-2} \frac{|\omega|\sqrt{1 + 2 \omega^2}}{(1 + \omega^2)^{1/4}} \left| \widehat{h_1}(\omega\sqrt{\omega^2 + 1})\right| \d \omega \right\|_{L^2_{\tau}} \\
&\lesssim \left\| \int _{-\infty}^\infty \lb \tau - \omega \sqrt{\omega^2 + 1}\rb^{-2} \lb \omega \sqrt{\omega^2 + 1} \rb^{1/4} \left| \widehat{h_1}(\omega\sqrt{\omega^2 + 1})\right| \frac{2 \omega^2 + 1}{\sqrt{\omega^2 + 1}}\d \omega \right\|_{L^2_{\tau}} \\
&= \left\| \int _{-\infty}^\infty \lb \tau - z\rb^{-2} \lb z \rb^{1/4} \left| \widehat{h_1}(z)\right| \d z \right\|_{L^2_{\tau}} \\
&\lesssim \int_{-\infty}^\infty \lb z \rb^{1/2} \left|\widehat{h_1}(z)\right|^2 \d z = \| \chi h_1\|_{H^{\frac14}(\R)}.
\end{align*}
The last line follows from an application of Young's inequality. 

The procedure for $B$ is exactly the same -- we drop the factor of $\omega$ and replace $h_1$ with $h_2$ in the integrals above to arrive at a bound of 
\[ \int_{-\infty} ^\infty \lb z \rb^{-3/2} \left|\widehat{h_2}(z)\right|^2 \d z = \| \chi h_2\|_{H^{-\frac14}(\R)}.\] 

It remains to obtain bounds on $A$ and $B$ in $X^{s,\frac12-}$ for general $s$. Notice that for any $s \in \N$, the derivative $\partial_x^s\bigl(\eta A\bigr)$ is 
\[ \eta(t) \int_{-\infty}^\infty \frac{i\omega\bigl(i\omega + \sqrt{\omega^2 + 1}\bigr)}{\sqrt{\omega^2 + 1}} e^{it \omega \sqrt{\omega^2 + 1}} \; \widehat{h_1}\Bigl( \omega \sqrt{\omega^2 + 1}\Bigr) f^{(s)}\Bigl(x \sqrt{\omega^2 + 1}\Bigr)  (\omega^2 +1 )^{s/2} \d \omega, \]
with a similar formula for $\partial_x^s\bigl(\eta B\bigr)$. Since $(\omega^2 + 1)^{s/2} \lesssim \lb \omega \sqrt{\omega^2 +1} \rb^{s/2}$, the desired result follows for $s \in \N$. By interpolation, we obtain the bound for any $s > 0$. 

For $s < 0$, let $\lb \partial \rb ^{-1/2}_x$ be the Fourier multiplier operator $\lb \xi \rb^{-1/2}$. Then $\lb \partial \rb^{-1/2}_x \bigl(\eta A\bigr)$ is equal to 
\[  \eta(t) \int_{-\infty}^\infty \frac{i\omega\bigl(i\omega + \sqrt{\omega^2 + 1}\bigr)}{\sqrt{\omega^2 + 1}} e^{it \omega \sqrt{\omega^2 + 1}} \; \widehat{h_1}\Bigl( \omega \sqrt{\omega^2 + 1}\Bigr) \lb \partial \rb^{-1/2}_x \Bigl[ f\Bigl(x \sqrt{\omega^2 + 1}\Bigr) \Bigr] \d \omega, \]
again with a similar statement for $B$. Now notice that 
\begin{align*}
\mathcal{F}_x \Biggl( \lb \partial \rb^{-1/2}_x \Bigl[ f\Bigl(x \sqrt{\omega^2 + 1}\Bigr) \Bigr] \Biggl) =  \mathcal{F}_x \Biggl( \Bigl[ \lb \partial \rb^{-1/2}_x f \Bigr] \Bigl(x \sqrt{\omega^2 + 1}\Bigr) \Biggl)(\xi) \; \frac{\lb \xi/ \sqrt{\omega^2 + 1} \rb^{1/2} }{\lb \xi \rb^{1/2}}. 
\end{align*}
Noting the $\lb \partial \rb_x^{-1/2}f$ is also a Schwarz function, we proceed just as in the case $s = 0$. In that situation, we moved the $L^2_\xi$ norm inside the integral and used the fact that $\| \frac{1}{(1 + \omega^2 + \xi^2)^{1/2+}}\|_{L^2_\xi} \approx (1 + \omega^2)^{-1/4-} \leq (1+\omega^2)^{-1/4}$. In this case, we use
\[ \left\| \frac{\lb \xi/ \sqrt{\omega^2 + 1} \rb^{1/2}}{\lb \xi \rb^{1/2} ( 1 + \omega^2 + \xi^2)^{1/2+}} \right\|_{L^2_\xi} = \left\| \frac{1}{(1 + \omega^2)^{1/4} (1 + \xi^2)^{1/4} ( 1 + \omega^2 + \xi^2)^{1/4+}} \right\|_{L^2_\xi} \lesssim \frac{1}{(1 + \omega^2)^{1/2}}.\]
This bound holds since 
\[ \int \frac{1}{(1 + \xi^2)^{1/2} (1 + \omega^2 + \xi^2)^{1/2+}} \d \xi \approx  \int \frac{1}{\lb  \xi\rb \lb |\omega| + |\xi| \rb^{1+}} \d \xi \lesssim \lb \omega \rb^{-1} \]
by Lemma \ref{calc_est}. 
Then the same argument we used previously yields the bound 
\[ \int_{-\infty} ^\infty \left|\widehat{h_1}(z)\right|^2 \d z.\]
We obtain a similar bound for $B$. Interpolating between the $s = -\frac12$ and $s=0$ estimates completes the proof. 

\subsection{Proof of Lemma \ref{cont}: Continuity of Linear Flow} \label{cont_prf}

Recall that $2 \pi \; W_0^t(0,0,h_1,h_2) = -A -B  + C +D$. We start with the claim $A,B  \in C_t^0H^s_x(\R \times \R)$. Note that   
\begin{align*}
A = \int_{-\infty}^\infty f\left(x\sqrt{\omega^2 + 1}\right) \mathcal{F}(L^t\phi_A)(\omega) \d \omega ,\\
B = \int_{-\infty}^\infty f\left(x\sqrt{\omega^2 + 1}\right) \mathcal{F}(L^t\phi_B)(\omega) \d \omega
\end{align*}
where 
\begin{align*}
\widehat{\phi_A} = \frac{i\omega\bigl(i\omega + \sqrt{1 + \omega^2}\bigr)}{\sqrt{\omega^2+1}} \widehat{h_1}\left(\omega \sqrt{\omega^2 + 1}\right), \quad \quad 
\widehat{\phi_B} = \frac{i\omega + \sqrt{1 + \omega^2}}{\sqrt{\omega^2+1}} \widehat{h_2}\left(\omega \sqrt{\omega^2 + 1}\right), 
\end{align*}
the function $f$ is given by $f(x) = e^{-x} \rho(x)$, and $L^t$ is the Fourier multiplier operator with multiplier $e^{i\omega\sqrt{\omega^2 + 1} t}$. Now 
\begin{align*}
\| \phi_A \|_{H^s_x}^2 
&= \int_{-\infty}^\infty \frac{\omega^2(2\omega^2 + 1)}{\omega^2 + 1} \bigl\lb \omega \bigr\rb^{2s} \; \left|\widehat{h_1}\Bigl(\omega \sqrt{\omega^2 + 1}\Bigr)\right|^2 \d \omega \\
&= \int_{-\infty}^\infty  \frac{\omega^2 \lb \omega \rb^{2s}\bigl\lb \omega \sqrt{\omega^2 + 1} \bigr\rb^{\frac{2s + 1}{2}}}{\sqrt{\omega^2 + 1} \bigl\lb \omega \sqrt{\omega^2 + 1} \bigr\rb^{\frac{2s + 1}{2}} }  \; \left| \widehat{h_1}\Bigl(\omega \sqrt{\omega^2 + 1}\Bigr) \right|^2 \frac{2\omega^2 + 1}{\sqrt{\omega^2 + 1}} \d \omega  \\
&\lesssim \int_{-\infty}^\infty \bigl\lb \omega \sqrt{\omega^2 + 1} \bigr\rb^{\frac{2s + 1}{2}} \; \left| \widehat{h_1}\Bigl(\omega \sqrt{\omega^2 + 1}\Bigr) \right|^2 \frac{2\omega^2 + 1}{\sqrt{\omega^2 + 1}} \d \omega \\
&= \int_{-\infty}^\infty \lb z \rb^{\frac{2s + 1}{2}} \; \left| \widehat{h_1}(z) \right|^2 \d z = \| \chi h_1 \|_{H^\frac{2s+1}{4}(\R)}^2,
\end{align*}
and similarly 
\[ \| \phi_B \|_{H^s_x}^2 \lesssim \| \chi h_2 \|_{H^\frac{2s-1}{4}(\R)}^2. \] 
Thus, using time continuity of the linear operator $L^t$, it suffices to show that the map 
\[ g \mapsto T(g) = \int_{-\infty}^\infty f\left(x \sqrt{\omega^2 + 1}\right) \widehat{g}(\omega) \d \omega \]
is bounded from $H^s$ to $H^s$. Consider first $s = 0$. 

Rewrite $Tg(x)$ as follows using the change of variables $z = x \sqrt{\omega^2 + 1}$: 
\begin{align*}
Tg(x) &= \int _{-\infty} ^\infty f\left(x\sqrt{\omega^2 + 1}\right) \widehat{g}(\omega) \d \omega \\
&= \int_x^{\operatorname{sgn}(x)\infty} f(z)\;  \Bigl[ \widehat{g}\Bigl( \sqrt{(z/x)^2 - 1} \Bigr) + \widehat{g}\Bigl( -\sqrt{(z/x)^2 - 1} \Bigr) \Bigr]  \frac{z/x^2}{\sqrt{(z/x)^2 - 1}} \d z.
\end{align*}

Then 
\begin{align*}
\|Tg\|_{L^2_x} &\lesssim \int |f(z)| \left\| \chi_{[0,z]}(x) \; \widehat{g}\left(\pm \sqrt{(z/x)^2 - 1} \right) \frac{z/x^2}{\sqrt{(z/x)^2 - 1}} \right\|_{L^2_x} \d z,
\end{align*}
and, expanding the $L^2_x$ norm, 
\[ \int_0^z \left| \widehat{g} \left(\pm  \sqrt{(z/x)^2 - 1} \right)\right|^2 \frac{z^2/x^4}{(z/x)^2 -1} \d x = \frac{1}{|z|} \int_0^{\pm\infty} \left| \widehat{g}(y) \right|^2 \frac{\sqrt{1 + y^2}}{y} \d y .\] 
On the region where $|y| \geq 1$, the right hand side above is bounded by $\frac{1}{|z|} \|g\|_{L^2}^2$. Since $|f(z)|/\sqrt{|z|}$ is in $L^1$, this yields the desired bound. For the case when $|y| \leq 1$, go back to the form $Tg(x) = \int  f\left(x\sqrt{\omega^2 + 1}\right) \widehat{g}(\omega) \d \omega $. The region $|y| \leq 1$ corresponds to $|\omega| \leq 1$ in this integral. So we consider the following norm, which we bound by applying Cauchy-Schwarz in $\omega$ and then using the change of variables $y = x \sqrt{1 + \omega^2}$ to replace the integration in $x$:
\begin{align*}
\left\| \int _{-1} ^1 f\left(x\sqrt{\omega^2 + 1}\right) \widehat{g}(\omega) \d \omega \right\|_{L^2_x}^2 &\lesssim \|g\|_{L^2}^2 \left\| \chi_{[0,1]}(\omega) f\left(x\sqrt{\omega^2 + 1}\right) \right\|_{L^2_{x,\omega}}^2 \\
&= \|g\|_{L^2}^2 \int_{-1}^1 \frac{1}{\sqrt{1 + \omega^2}}\int f^2(y) \d y \d \omega \\
&\lesssim \|g\|_{L^2}^2.
\end{align*}
This completes the proof that $A, B \in C^0_tH^s_x$ for $s = 0$.  

For $s > 0 $, notice that for any $s \in \N$, we have 
\begin{align*}
\partial_x^s Tg(x) = \int_0^\infty f^{(s)}\left(x\sqrt{\omega^2 + 1}\right) (\omega^2 + 1)^{s/2}\;  \widehat{g}(\omega) \d \omega. 
\end{align*}
This and interpolation imply the desired bounds for $A$ and $B$ in $H^s_x$ for positive $s$. 

Also, if we choose $\rho$ such that $\int f \d x = 0$ so that $\partial_x^{-1}f$ is a Schwarz function, then we have
\[  \partial_x^{-1} Tg(x) = \int_0 ^\infty  \partial_x^{-1} f\left(x\sqrt{\omega^2 + 1}\right) (\omega^2 + 1)^{-1/2}\;  \widehat{g}(\omega) \d \omega. \] 
Combining this with the $s=0$ result and interpolation, we obtain the bound for $s \geq -1$. 

Next, recall that 
\begin{align*}
C &= L^t\phi_C(x) \quad \text{where} \quad \widehat{\phi_C}(\omega) = \bigl(i \omega + \sqrt{\omega^2 + 1}\bigr) \;\widehat{h_1}\left( \omega \sqrt{\omega^2 + 1}\right), \\
D &= L^t\phi_D(x) \quad \text{where} \quad \widehat{\phi_D}(\omega) =  \frac{i \omega + \sqrt{\omega^2 + 1}}{\sqrt{\omega^2 + 1}} \;\widehat{h_2}\left( \omega \sqrt{\omega^2 + 1}\right).
\end{align*}
The $C^0_tH^s_x$ bounds for these terms follow from the continuity of the linear operator $F^t$ and the bounds for $\phi_C$ and $\phi_D$ which were proved in Lemma \ref{linear_half}.

It remains to prove that $\eta W_0^t(0,0,h_1,h_2)$ is in $C^0_xH^{\frac{2s+1}{4}}_t$. Recall the form of $C$ and $D$ as linear flows and apply Lemma \ref{Kato} to obtain the desired bound for these terms. For $A$ and $B$, write
\begin{align*}
A(x,t) = \int \mathcal{F}_\omega\left(  f\bigl(x\sqrt{\omega^2 + 1}\bigr)\right) \!(y) \;   L^t\phi_A(y) \d y ,\\
B(x,t) = \int \mathcal{F}_\omega \left( f\bigl(x\sqrt{\omega^2 + 1}\bigr) \right)\!(y) \;   L^t\phi_B(y) \d y
\end{align*}
where 
$\phi_A$ and $\phi_B$ are defined as before. Then $A$ is equal to 
\[ \int \!\frac1x \mathcal{F}_z\Bigl(f(\operatorname{sgn}(x) \sqrt{x^2  +z^2})\Bigr)\left(\frac{y}x\right) L^t\phi_A(y) \d y =  \int \mathcal{F}_z\Bigl(f(\operatorname{sgn}(x) \sqrt{x^2  +z^2})\Bigr)(y) L^t\phi_A(xy) \d y, \]
with a parallel statement for $B$. By Kato smoothing, Lemma \ref{Kato}, it suffices to show that the function $\mathcal{F}_z\left( f\bigl(\operatorname{sgn}(x) \sqrt{z^2 + x^2} \bigr)\right)(y)$ is in $L^\infty_xL^1_y$. It is enough to show that $f\bigl(\operatorname{sgn}(x) \sqrt{z^2 + x^2} \bigr)$ is in $L^\infty_xH^1_z$ since
\[ \int |\widehat{k}(y)| \d y = \int \lb y \rb |\widehat{k}(y)| \lb y \rb^{-1} \d y \lesssim \|k \|_{H^1} \| \lb \cdot \rb^{-1} \|_{L^2}.\] 
To this end, we consider the $L^2_z$ and the $\dot{H}^1_z$ norms separately. For the $L^2_z$ norm, split the integral into two regions, one where $|z|$ is small, and its complement:
\begin{align*}
&\;\; \int \left|f\Bigl(\operatorname{sgn}(x) \sqrt{z^2 + x^2} \Bigr)\right|^2 \d z \\
&= \int_{|z| \leq 1} \left|f\Bigl(\operatorname{sgn}(x) \sqrt{z^2 + x^2} \Bigr)\right|^2 \d z +  \int_{|z| > 1} \left|f\Bigl(\operatorname{sgn}(x) \sqrt{z^2 + x^2} \Bigr)\right|^2 \d z.
\end{align*}
The first term is bounded since $f$ is bounded. Set $y = \operatorname{sgn}(x)\sqrt{z^2 + x^2}$ in the second integral to obtain
\begin{align*}
\int_{|y^2-x^2| > 1} \left|f(y)\right|^2 \frac{y}{\sqrt{y^2 - x^2}} \d y,
\end{align*}
which is bounded since $f$ is a Schwarz function. The same argument serves to bound the derivative since \[
\frac{\d}{\d z} \left[ f\Bigl(\operatorname{sgn}(x) \sqrt{z^2 + x^2} \Bigr)\right] = f'\Bigl(\operatorname{sgn}(x) \sqrt{z^2 + x^2} \Bigr) \frac{\operatorname{sgn}(x) z}{\sqrt{z^2 + x^2}} \quad \text{ and } \quad  \frac{ |z|}{\sqrt{z^2 + x^2}} \leq 1. \] 

\subsection{Proof of Lemma \ref{bilin}: Bilinear $X^{s,b}$ Estimate} \label{bilin_prf}

By duality, it suffices to show that 
\begin{equation} \label{eq:dualForm}
\left| \iint \mathcal{M}(uv)_{xx} \overline{\phi} \d x \d t \right| \lesssim \| u\|_{X^{s,b}}\| v\|_{X^{s,b}} \| \phi\|_{X^{-(s+a), b}} 
\end{equation}
for any $\phi \in X^{-(s+a), b}$. 
The left-hand side of \eqref{eq:dualForm} is equal to 
\begin{align*}
\left| \iint \frac{\xi^2 \; \widehat{uv}(\xi,\tau) }{\sqrt{\xi^2 + \xi^4}} \; \overline{\widehat{\phi}(\xi,\tau)} \d\xi \d\tau \right| = \left| \iiiint \frac{\xi^2 \; \widehat{u}(\xi_1, \tau_1) \widehat{v}(\xi - \xi_1,\tau-\tau_1) }{\sqrt{\xi^2 + \xi^4}} \; \overline{\widehat{\phi}(\xi,\tau)} \d \xi_1 \d \tau_1 \d\xi \d\tau \right|. 
\end{align*}
Now define
\begin{align*} p(\xi, \tau) = \lb \xi \rb^s \bigl\lb |\tau| - \sqrt{\xi^2 + \xi^4}& \bigr\rb^b \; \widehat{u}(\xi,\tau), 
\qquad \qquad  q(\xi, \tau) = \lb \xi \rb^s \bigl\lb |\tau| - \sqrt{\xi^2 + \xi^4} \bigr\rb^b \; \widehat{v}(\xi,\tau), \\
 r(\xi, \tau) &= \lb \xi \rb^{-(s+a)} \bigl\lb |\tau| - \sqrt{\xi^2 + \xi^4} \bigr\rb^b \; \overline{\widehat{\phi}(\xi,\tau)}. 
\end{align*}
The desired bound \eqref{eq:dualForm} is equivalent to showing that 
\begin{align*}
\left| \iiiint M(\xi, \xi_1, \tau, \tau_1) p(\xi_1, \tau_1) q(\xi - \xi_1, \tau - \tau_1) r(\xi, \tau) \d \xi_1 \d \tau_1 \d \xi  \d \tau \right| \lesssim \|p\|_{L^2_{\xi,\tau}} \|q\|_{L^2_{\xi,\tau}} \|r\|_{L^2_{\xi,\tau}},
\end{align*}
where the multiplier ${M}$ is 
\[ {M} = \frac{\xi^2 \lb \xi \rb^{s + a} \lb \xi_1 \rb ^{-s} \lb \xi - \xi_1 \rb^{-s} }{\sqrt{ \xi^2 + \xi^4} \lb |\tau| - \xi^2 \rb ^{b} \lb |\tau_1 | - \xi_1^2 \rb^b \lb | \tau - \tau_1 | - (\xi - \xi_1)^2 \rb^b }. \] 
There are six possibilities for the signs of $\tau$, $\tau_1$, and $\tau - \tau_1$:
\begin{enumerate}[label=(\alph*) ]
\item \label{1} $\tau_1 \geq 0$, $\tau - \tau_1 \geq 0$,
\item \label{2} $\tau_1 \geq 0$, $\tau - \tau_1 \leq 0$, and $\tau \geq 0$,
\item \label{3} $\tau_1 \geq 0$, $\tau - \tau_1 \leq 0$, and $\tau \leq 0$,
\item \label{4} $\tau_1 \leq 0$, $\tau - \tau_1 \leq 0$,
\item \label{5} $\tau_1 \leq 0$, $\tau - \tau_1 \geq 0$, and $\tau \leq 0$,
\item \label{6} $\tau_1 \leq 0$, $\tau - \tau_1 \geq 0$, and $\tau \geq 0$.
\end{enumerate}
Since $L^2$ norms are invariant under reflections, we can use the substitution $(\tau, \tau_1) \mapsto -(\tau, \tau_1)$ to reduce \ref{4}, \ref{5}, and \ref{6} to \ref{1}, \ref{2}, and \ref{3} respectively. 

Consider first \ref{1}. By Cauchy-Schwarz in the $\xi$-$\tau$ integral, it suffices to show that
\[ \left\| \iint M(\xi, \xi_1, \tau, \tau_1) p(\xi_1, \tau_1) q(\xi - \xi_1, \tau - \tau_1) \d \xi_1 \d \tau_1 \right\|_{L^2_{\xi, \tau}} \lesssim \|p\|_{L^2_{\xi,\tau}} \|q\|_{L^2_{\xi,\tau}}.\] 
Using Cauchy-Schwarz and Young's inequalities, the left-hand side of this is bounded by 
\begin{align*}
&\left\| \|M\|_{L^2_{\xi_1, \tau_1}} \|f(\xi_1, \tau_1)g(\xi - \xi_1, \tau -\tau_1) \|_{L^2_{\xi_1,\tau_1}} \right\|_{L^2_{\xi,\tau}} \\
&\qquad \lesssim \; \Bigl( \sup_{\xi, \tau} \|M\|_{L^2_{\xi_1, \tau_1}} \Bigr) \|f(\xi_1, \tau_1)g(\xi - \xi_1, \tau -\tau_1) \|_{L^2_{\xi_1,\tau_1, \xi,\tau}} \\
&\qquad = \; \Bigl(\sup_{\xi, \tau} \|M\|_{L^2_{\xi_1, \tau_1}} \Bigr) \| f^2 * g^2 \|_{L^1_{\xi, \tau}}^{1/2} \\
&\qquad \lesssim \; \Bigl(\sup_{\xi, \tau} \|M\|_{L^2_{\xi_1, \tau_1}} \Bigr) \|f\|_{L^2_{\xi, \tau}} \|g\|_{L^2_{\xi,\tau}}.
\end{align*}

Thus, in Case \ref{1}, it suffices to show that 
\[ \sup_{\xi, \tau}  \iint \frac{\xi^4 \lb \xi \rb^{2s + 2a} \lb \xi_1 \rb ^{-2s} \lb \xi - \xi_1 \rb^{-2s}}{(\xi^2 + \xi^4) \lb \tau - \xi^2 \rb ^{2b} \lb \tau_1  - \xi_1^2 \rb^{2b} \lb \tau - \tau_1  - (\xi - \xi_1)^2 \rb^{2b} } \d \xi_1 \d \tau_1  \] 
is finite. Using the fact that $\lb a \rb \lb b \rb \gtrsim \lb a+b \rb$, we can eliminate the $\tau$ dependence to obtain
\[ \sup_{\xi}  \iint \frac{\xi^4 \lb \xi \rb^{2s + 2a} \lb \xi_1 \rb ^{-2s} \lb \xi - \xi_1 \rb^{-2s}}{(\xi^2 + \xi^4) \lb \tau_1  - \xi_1^2 \rb^{2b} \lb \tau_1 - 2\xi\xi_1 + \xi_1^2 \rb^{2b} } \d \xi_1 \d \tau_1  \] 
Applying Lemma \ref{calc_est} in $\tau_1$ and observing that $\xi^4/(\xi^2 + \xi^4) < 1$, we are reduced to bounding 
\[ \sup_{\xi} \int \frac{\lb \xi \rb^{2s + 2a} \lb \xi_1 \rb^{-2s} \lb \xi - \xi_1 \rb^{-2s}}{\lb \xi_1(\xi_1 - \xi) \rb^{1-}} \d \xi_1. \]
We consider several cases.\\
\textbf{Case 1.} $\lb \xi \rb^{s} \lesssim \lb \xi_1 \rb^s \lb \xi - \xi_1 \rb^s$. \\
This case reduces to bounding 
\[ \sup_\xi \; \lb \xi \rb^{2a} \int \frac{1}{\lb \xi_1 (\xi_1 - \xi) \rb ^{1-}} \d \xi_1. \]
Let 
\[ x = \xi_1(\xi-\xi_1) \quad \Rightarrow \quad 2 \xi_1 = \xi \pm \sqrt{\xi^2 + 4x} \quad \text{and} \quad \d x = \pm \sqrt{\xi^2 + 4x} \d \xi_1. \]
Then the supremum above is bounded by 
\[ \sup_\xi \; \lb \xi \rb^{2a} \int \frac{1}{\lb x\ \rb^{1-} \sqrt{\xi^2 + 4x}} \d x. \]
By \cite[Lemma 6.3]{ET}, this is bounded by 
\[ \sup_\xi \; \lb \xi \rb^{2a - 1 +}, \]
which is finite as long as $a < \frac12$. \\
\textbf{Case 2.} $\lb \xi \rb \ll \lb \xi_1 \rb \lb \xi - \xi_1 \rb $ and $s < 0$. \\
\textbf{Case 2a.} $|\xi|   \lesssim 1$. \\
In this case, we must control $\int \lb \xi_1 \rb ^{-4s - 2 + } \d \xi_1$, which is possible when $s > -\frac14$. \\
\textbf{Case 2b.} $|\xi_1| \gg |\xi| \gtrsim 1$. \\
In this case, we arrive at 
\[ \sup_{\xi} \; \lb \xi \rb^{2s + 2a} \int_{|\xi_1| \gg |\xi|} \lb \xi_1 \rb^{-4s - 2 +} \d \xi_1 \; \lesssim \; \sup_{\xi} \; \lb \xi \rb ^{2a -2s - 1 +}.\]
This is finite when $a < s + \frac12$. \\
\textbf{Case 2c.} $|\xi| \gg |\xi_1| \gg 1$. \\
In this case, we arrive at
\[ \sup_{\xi} \; \lb \xi \rb^{2a - 1 +}\int_{|\xi_1| \ll |\xi|} \lb \xi_1 \rb^{-2s-1 + } \d \xi_1. \]
Since $s < 0$, this converges if $a < s + \frac12$. \\
\textbf{Case 2d.} $|\xi_1| \approx  |\xi| $. \\
This is only possible if $|\xi - \xi_1| \gg 1$ and $|\xi_1| \gg 1$. Thus we need to bound 
\[ \sup_\xi \; \lb \xi \rb^{2a} \int \lb \xi - \xi_1 \rb^{-1-2s +} \lb \xi_1 \rb^{-1+} \d \xi_1. \]
Using Lemma \ref{calc_est}, we see that this can be bounded for $s > -\frac12$ as long as $a < s + \frac12$. 
This completes the proof for the combination of $\tau$ signs described in \ref{1}. 

For the combination of signs described in \ref{2}, we follow the same procedure of estimating using Cauchy-Schwarz and Young's inequalites, but exchange the role of $(\xi, \tau)$ and $(\xi_1, \tau_1)$. It then suffices to control
\begin{equation} \label{eq:secondCase}
\sup_{\xi_1} \int \frac{\lb \xi \rb^{2s + 2a} \lb \xi_1 \rb^{-2s} \lb \xi - \xi_1 \rb^{-2s}}{\lb \xi(\xi_1 - \xi) \rb^{1-}} \d \xi. 
\end{equation}
Case \ref{3} can be reduced to the same estimate by performing the change of variables $(\xi_1, \tau_1) \mapsto (\xi - \xi_1, \tau - \tau_1)$ and then carrying out the same series of estimates. To bound this supremum \eqref{eq:secondCase}, we consider similar cases. \\
\textbf{Case 1'.} $\lb \xi \rb^{s} \lesssim \lb \xi_1 \rb^s \lb \xi - \xi_1 \rb^s$. \\
The procedure here is precisely the same as in \textbf{Case  1}. \\
\textbf{Case 2'.} $\lb \xi \rb \ll \lb \xi_1 \rb \lb \xi - \xi_1 \rb $ and $s < 0$. \\
\textbf{Case 2a'.} $|\xi_1|   \lesssim 1$. \\
In this case, we must control $\int \lb \xi_1 \rb ^{2a - 2 + } \d \xi_1$, which is possible when $a < \frac12$. \\
\textbf{Case 2b'.} $|\xi_1| \gg |\xi| \gtrsim 1$. \\
In this case, we arrive at 
\[ \sup_{\xi_1} \; \lb \xi_1 \rb^{-4s -1+} \int_{|\xi| \ll |\xi_1|} \lb \xi \rb^{2s + 2a -1 +} \d \xi \; \lesssim \; \sup_{\xi_1} \; \lb \xi_1 \rb ^{-4s - 1 + \max\{ 2s + 2a, 0\}+},\]
which is finite if $a < s + \frac12$ and $s > -\frac14$. \\
\textbf{Case 2c'.} $|\xi| \gg |\xi_1| \gg 1$. \\
In this case, we arrive at
\[ \sup_{\xi_1} \; \lb \xi_1 \rb^{-2s}\int_{|\xi| \gg |\xi_1|} \lb \xi \rb^{2a-2+ } \d \xi \; \lesssim \; \sup_{\xi_1} \; \lb \xi_1 \rb^{-2s+ 2a - 1 +}, \]
which is finite if $a < s + \frac12$. \\
\textbf{Case 2d'.} $|\xi_1| \approx  |\xi| $. \\
This case only arises if $|\xi - \xi_1| \gg 1$ and $|\xi_1| \gg 1$. Thus we need to bound 
\[ \sup_{\xi_1} \; \lb \xi_1 \rb^{-2s} \int \lb \xi - \xi_1 \rb^{-2s -1  +} \lb \xi \rb^{2s+2a-1+} \d \xi  . \]
Using Lemma \ref{calc_est}, we see that this can be bounded for $ -\frac12< s < 0$ and $a < \frac12$ by 
\[ \sup_{\xi_1} \; \lb \xi_1 \rb^{-2s + \max\{ -2s, 2s + 2a\} - 1 +}.\] This is finite for $a < \frac12$ and $s > -\frac14$.  
This completes the proof.

\subsection{Proof of Lemma \ref{Duhamel_Kato}: Kato Smoothing for Duhamel Term} \label{Duhamel_Kato_prf} 

Again, it suffices to consider evaluation at $x = 0$ since a spatial translation of $G$ does not affect the magnitude of $\widehat{\mathcal{M}(G)}$. At $x= 0$, we have 
\[ \int_0^t W_{R,2}^{t-t'} G \d t' = \frac{1}{2i} \int \int_0^t \frac{ e^{i(t-t') \sqrt{\xi^2 + \xi^4}} - e^{-i(t-t') \sqrt{\xi^2 + \xi^4}}}{\sqrt{\xi^2 + \xi^4}} \mathcal{F}_x(G)(\xi, t') \d t' \d \xi.\] 
Also, note that 
\[ \mathcal{F}_x(G)(\xi, t') = \int e^{i\tau t'} \widehat{G}(\xi, \tau) \d \tau \]
and 
\[ \int_0^t e^{it'(\tau \mp \sqrt{\xi^2 + \xi^4})} \d t' = \frac{e^{it(\tau \mp \sqrt{\xi^2 + \xi^4})}- 1}{i(\tau \mp \sqrt{\xi^2 + \xi^4})}. \]
Thus we wish to bound 
\begin{align*}
\iint \frac{e^{it\tau} - e^{\pm it \sqrt{\xi^2 + \xi^4}}}{\sqrt{\xi^2 + \xi^4}\left(\tau \mp \sqrt{\xi^2 + \xi^4}\right)} \widehat{G}(\xi, \tau) \d \xi \d \tau. 
\end{align*}
Let $\Psi$ be a smooth cut-off function such that $\Psi = 1$ on $[-1,1]$ and $\Psi = 0$ outside $[-2,2]$. Let $\Psi^C = 1 - \Psi$. Then write 
\begin{align*}
2\eta(t) \int_0^t F_2^{t-t'} G \d t' &= \eta(t) \iint \frac{\Bigl(e^{it\tau} - e^{\pm it \sqrt{\xi^2 + \xi^4}}\Bigr) \; \Psi\Big(\tau \mp \sqrt{\xi^2 + \xi^4}\Big)}{\sqrt{\xi^2 + \xi^4}\left(\tau \mp \sqrt{\xi^2 + \xi^4}\right)} \widehat{G}(\xi, \tau)  \d \xi \d \tau \\
&+ \eta(t) \iint \frac{e^{it\tau}\; \Psi^C\Big(\tau \mp \sqrt{\xi^2 + \xi^4}\Big)}{\sqrt{\xi^2 + \xi^4}\left(\tau \mp \sqrt{\xi^2 + \xi^4}\right)} \widehat{G}(\xi, \tau)  \d \xi \d \tau \\
&- \eta(t) \iint \frac{ e^{\pm it \sqrt{\xi^2 + \xi^4}}\Psi^C\Big(\tau \mp \sqrt{\xi^2 + \xi^4}\Big)}{\sqrt{\xi^2 + \xi^4}\left(\tau \mp \sqrt{\xi^2 + \xi^4}\right)} \widehat{G}(\xi, \tau) \d \xi \d \tau \\
&= \I + \II - \III.
\end{align*}
By Taylor expanding, we have 
\begin{align*}
\frac{e^{it\tau} - e^{\pm it \sqrt{\xi^2 + \xi^4}}}{\left(\tau \mp \sqrt{\xi^2 + \xi^4}\right)} &= -e^{it\tau} \sum_{k=1}^\infty \frac{(-it)^k}{k!} (\tau \mp \sqrt{\xi^2 + \xi^4})^{k-1}. 
\end{align*}
Therefore $\| \I\|_{H^{\frac{2s + 1}{4}}_t}$ is bounded by 
\begin{align*}
& \sum_{k=1}^\infty \frac{\| \eta(t) t^k\|_{H^1}}{k!} \left\| \iint e^{it\tau} \left(\tau \mp \sqrt{\xi^2 + \xi^4}\right)^{k-1} \Psi\Big(\tau \mp \sqrt{\xi^2 + \xi^4}\Big) \frac{\widehat{G}(\xi, \tau)}{\sqrt{\xi^2 + \xi^4}} \d \xi \d \tau \right\|_{H^{\frac{2s +1}{4}}_t} \\
&\lesssim \sum_{k=1}^\infty \frac{1}{(k-1)!}\left\| \lb \tau \rb^{\frac{2s + 1}{4}} \int \left(\tau \mp \sqrt{\xi^2 + \xi^4} \right)^{k-1} \Psi\Big(\tau \mp \sqrt{\xi^2 + \xi^4}\Big) \frac{ \widehat{G}(\xi, \tau)}{\sqrt{\xi^2 + \xi^4}} \d \xi \right\|_{L^2_\tau} \\
&\lesssim \left\| \lb \tau \rb^{\frac{2s + 1}{4}} \int \Psi\Big(\tau \mp \sqrt{\xi^2 + \xi^4}\Big) \frac{ \widehat{G}(\xi, \tau)}{\sqrt{\xi^2 + \xi^4}} \d \xi \right\|_{L^2_\tau}.
\end{align*}
Using the Cauchy-Schwarz inequality in $\tau$, this can be bounded by
\begin{align*}
&\;\; \left[ \int \lb \tau \rb^\frac{2s + 1}{2} \left( \int_{\left|\tau \mp \sqrt{\xi^2 + \xi^4}\right| < 1} \lb \xi \rb^{-2s} \d \xi \right) \left( \int_{\left|\tau \mp \sqrt{\xi^2 + \xi^4}\right| < 1} \lb \xi \rb^{2s} \frac{| \widehat{G}(\xi, \tau)|^2}{\xi^2 + \xi^4} \d \xi \right) \d \tau \right]^{1/2} \\
&\lesssim \; \sup_\tau \left( \lb \tau \rb^\frac{2s+1}{2} \int_{\left|\tau \mp \sqrt{\xi^2 + \xi^4}\right| < 1} \lb \xi \rb^{-2s} \d \xi\right)^{1/2}\|\mathcal{M} (G)\|_{X^{s,-b}}   \\
&\lesssim \|\mathcal{M} (G) \|_{X^{s,-b}}.
\end{align*}
The first inequality holds since on the region of interest in $\I$, we have
\[ 1 \approx \frac{1}{\bigl\lb \tau \mp \sqrt{\xi^2 + \xi^4} \bigr\rb^{b}} \leq \frac{1}{\bigl\lb |\tau| - \sqrt{\xi^2 + \xi^4} \bigr\rb^{b}}. \]
The supremum bound holds since 
\begin{align*}
\lb \tau \rb^\frac{2s+1}{2} \int_{\left|\tau \mp \sqrt{\xi^2 + \xi^4}\right| < 1} \lb \xi \rb^{-2s} \d \xi \;\lesssim \;
\begin{cases} 
1 & \text{ if } |\tau| \lesssim 1 \;\; \\
\lb \tau \rb^\frac{2s+1}{2} \int_{\left|\tau \mp \sqrt{|z| + z^2}\right| < 1} \lb z \rb^{-s - 1/2} \d z &\text{ if } |\tau| \gg 1.
\end{cases}   
\end{align*}
The latter bound comes from changing variables $\xi^2 \mapsto z$. The right-hand side is finite since the integrand is of order $|\tau|^{-s-1/2}$ over an interval of length $\approx 1$.  

Next consider $\III$. When $|\xi| \leq 1$, we have, using $b < \frac12$, the bound 
\begin{align*}
&\left\| \eta(t) \iint \frac{ e^{\pm it \sqrt{\xi^2 + \xi^4}}\Psi^C\Big(\tau \mp \sqrt{\xi^2 + \xi^4}\Big)}{\sqrt{\xi^2 + \xi^4}\left(\tau \mp \sqrt{\xi^2 + \xi^4}\right)} \widehat{G}(\xi, \tau) \d \xi \d \tau \right\|_{H^{\frac{2s+1}{4}}_t} \\
&\quad \lesssim \iint_{|\xi|\leq 1} \frac{\|\eta(t) e^{\pm it \sqrt{\xi^2 + \xi^4}}\|_{H^\frac{2s+1}{4}}}{\left| \tau \mp \sqrt{\xi^2 + \xi^4}\right| } \Psi^C\Big(\tau \mp \sqrt{\xi^2 + \xi^4}\Big) \left| \widehat{\mathcal{M}(G)}(\xi,\tau)\right| \d \xi \d \tau \\
&\quad \lesssim \iint \frac{\chi_{[-1,1]}(\xi)}{\bigl\lb \tau \mp \sqrt{\xi^2 + \xi^4}\bigr\rb } \left| \widehat{\mathcal{M}(G)}(\xi,\tau)\right| \d \xi \d \tau \\
&\quad \lesssim \|\mathcal{M}(G)\|_{X^{s,-b}} \Biggl\| \frac{\chi_{[-1,1]}(\xi)}{\bigl\lb \tau \mp \sqrt{\xi^2 + \xi^4}\bigr\rb^{1-b} } \Biggr\|_{L^2_{\tau, \xi}} \; \lesssim \;  \|\mathcal{M}(G)\|_{X^{s,-b}}. 
\end{align*}

To control the part of $\III$ where $|\xi| \geq 1$, change variables in the $\xi$ integral by setting $z = \pm \sqrt{\xi^2 + \xi^4} \approx \pm \xi^2$. Then, noticing that the integral is an inverse Fourier transform,  we obtain the bound
\begin{align*}
&\left\| \eta(t) \iint_{|\xi| \geq 1} \frac{ e^{\pm it \sqrt{\xi^2 + \xi^4}}\Psi^C\Big(\tau \mp \sqrt{\xi^2 + \xi^4}\Big)}{\sqrt{\xi^2 + \xi^4}\left(\tau \mp \sqrt{\xi^2 + \xi^4}\right)} \widehat{G}(\xi, \tau) \d \xi \d \tau  \right\|_{H^{\frac{2s+1}{4}}_t} \\
&\quad \lesssim \left\| \lb z \rb ^\frac{2s + 1}{4} \int \frac{| \widehat{G}(\xi(z), \tau) |}{\lb \tau - z\rb \bigl|4\xi(z)^3 + 2 \xi(z)\bigr|}   \d \tau \right\|_{L^2_{|z| \geq 2}} \\
&\quad \lesssim  \left\| \lb z \rb ^\frac{2s - 1}{4} \int \frac{ |\widehat{G}(\xi(z), \tau) |}{\lb \tau - z\rb \sqrt{\xi(z)^2 + \xi(z)^4}}   \d \tau \right\|_{L^2_{|z| \geq 2}}. 
\end{align*}
By Cauchy-Schwarz in the $\tau$ integral, using the fact that $b < \frac12$, this is bounded by 
\[ \left\| \lb z \rb ^\frac{2s - 1}{4} \frac{| \widehat{G}(\xi(z), \tau) |}{\lb \tau - z\rb^{b} (\xi(z)^2 + \xi(z)^4)^{1/2}}  \right\|_{L^2_{|z| \geq 2}L^2_\tau}. \] 
Changing variables back to $\xi$, this is bounded by $\|\mathcal{M}(G)\|_{X^{s,-b}}$, as desired. It remains to bound $\II$.

For $\II$, let $R$ denote the set $\{ |\tau| \gg |\xi|^2 \} \cup \{ |\xi| \lesssim 1\}$ and notice that 
\[\lb \tau \rb \; \lesssim \; \chi_R(\xi, \tau) \lb \tau - |\xi|^2 \rb + |\xi|^2 \]
and $(2s+1)/4 \geq 0$, so we have the bounds 
\begin{align*}
\| \II &\|_{H^\frac{2s+1}{4}_t} \lesssim \Bigl\| \lb \tau \rb^{\frac{2s+1}4} \int \frac{1}{\lb \tau \mp \sqrt{\xi^2 + \xi^4} \rb} \frac{| \widehat{G}(\xi,\tau)|}{\sqrt{\xi^2 +\xi^4}} \d \xi \Bigr\|_{L^2_\tau} \\
&\lesssim \Bigl\| \int \chi_R(\xi, \tau) \lb |\tau|  - \xi^2  \rb^\frac{2s-3}{4} \frac{| \widehat{G}(\xi,\tau)|}{\sqrt{\xi^2 +\xi^4}} \d \xi \Bigr\|_{L^2_\tau} 
+ \Bigl\| \int \frac{|\xi|^{s + 1/2}}{\bigl\lb \tau \mp \sqrt{\xi^2 + \xi^4} \bigr\rb} \frac{| \widehat{G}(\xi,\tau)|}{\sqrt{\xi^2 +\xi^4}} \d \xi \Bigr\|_{L^2_\tau} \end{align*}
The second term on the last line can be bounded by $\|\mathcal{M}(G)\|_{X^{s,b}}$ using Cauchy-Schwarz in $\xi$ provided that 
\[ \sup_\tau \int \frac{|\xi|}{\bigl\lb \tau \mp \sqrt{\xi^2 + \xi^4} \bigr\rb^{2-2b}} \d \xi \; < \; \infty, \]
which holds since $b < 1/2$ and 
\begin{equation*}
\sup_\tau \int \frac{|\xi|}{\bigl\lb \tau \mp \sqrt{\xi^2 + \xi^4} \bigr\rb^{2-2b}} \d \xi \;\lesssim\; \sup_\tau \int \frac{|\xi|}{\bigl\lb |\tau| - \xi^2 \bigr\rb^{2-2b}} \d \xi  \;\approx\; \sup_\tau \int \frac{1}{\bigl\lb |\tau| - z \bigr\rb^{2-2b}} \d z. 
\end{equation*}

For $s \leq \frac12$, we go back to 
\[ \Bigl\| \lb \tau \rb^{\frac{2s+1}4} \int \frac{1}{\bigl\lb \tau \mp \sqrt{\xi^2 + \xi^4} \bigr\rb} \frac{| \widehat{G}(\xi,\tau)|}{\sqrt{\xi^2 +\xi^4}} \d \xi \Bigr\|_{L^2_\tau} \]
and estimate using Cauchy-Schwarz in $\xi$ to obtain the bound
\[ \left[ \int \lb \tau \rb^{\frac{2s+1}{2}} \left( \int \frac{1}{\bigl\lb \tau \mp \sqrt{\xi^2 + \xi^4} \bigr\rb^{2-2b} \lb \xi \rb^{2s}} \d \xi \right) \left( \int \frac{\lb \xi \rb^{2s}}{\bigl\lb \tau \mp \sqrt{\xi^2 + \xi^4} \bigr\rb^{2b}} \frac{| \widehat{G}(\xi,\tau)|^2}{\xi^2 +\xi^4} \d \xi \right) \d \tau \right]^{1/2}, \]
which can be bounded by $\|\mathcal{M}(G)\|_{X^{s,-b}}$ as long as 
\[ \sup_\tau \; \lb \tau \rb^{\frac{2s+1}{2}} \int \frac{1}{\bigl\lb \tau \mp \sqrt{\xi^2 + \xi^4} \bigr\rb^{2-2b} \lb \xi \rb^{2s}} \d \xi\] 
is finite. To see that this holds for $s \leq \frac12$, recall that  $\lb \tau \pm \sqrt{\xi^2 + \xi^4} \rb \approx \lb \tau \pm \xi^2 \rb$. Consider $|\xi| \ll 1$ first. In this case, change variables in the integral by $\xi^2 \mapsto z$. Then apply Lemma \ref{calc_est} to bound the integral as follows:
\[  \sup_\tau \; \lb \tau \rb^{\frac{2s+1}{2}} \int \frac{1}{\lb \tau \mp z \rb^{2-2b} \lb z \rb^{s + \frac12}} \d z \lesssim \;  \sup_\tau \; \lb \tau \rb^{s + \frac12 - \min\{2-2b, s + \frac12\}} \; < \; \infty, \]
assuming that $b < \frac12$ and $ -\frac12 \leq s \leq \frac12$. Similarly, if $|\xi| \gtrsim 1$, again change variables by setting $z = \xi^2$ and apply Lemma \ref{calc_est} to obtain $\lb \tau \rb ^{s + \frac12 - 2 + 2b}$. This is finite for $b \leq \frac12$ and $s \leq \frac12$.

For the estimate on the derivative term, the procedure is similar. We break the Duhamel integral down into three pieces $\widetilde{\I} + \widetilde{\II} - \widetilde{\III}$:
\begin{align*}
\eta(t) \int_0^t W_{R,2}^{t-t'} G_{x} \d t' &= \eta(t) \iint \frac{i\xi \Bigl(e^{it\tau} - e^{\pm it \sqrt{\xi^2 + \xi^4}}\Bigr) \; \Psi\Big(\tau \mp \sqrt{\xi^2 + \xi^4}\Big)}{\sqrt{\xi^2 + \xi^4}\left(\tau \mp \sqrt{\xi^2 + \xi^4}\right)} \widehat{G}(\xi, \tau)  \d \xi \d \tau \\
&+ \eta(t) \iint \frac{i\xi e^{it\tau}\; \Psi^C\Big(\tau \mp \sqrt{\xi^2 + \xi^4}\Big)}{\sqrt{\xi^2 + \xi^4}\left(\tau \mp \sqrt{\xi^2 + \xi^4}\right)} \widehat{G}(\xi, \tau)  \d \xi \d \tau \\
&- \eta(t) \iint \frac{i\xi e^{\pm it \sqrt{\xi^2 + \xi^4}}\Psi^C\Big(\tau \mp \sqrt{\xi^2 + \xi^4}\Big)}{\sqrt{\xi^2 + \xi^4}\left(\tau \mp \sqrt{\xi^2 + \xi^4}\right)} \widehat{G}(\xi, \tau) \d \xi \d \tau \\
&= \widetilde{\I} + \widetilde{\II} - \widetilde{\III}.
\end{align*}
The only difference from the previous case is that each term now has a factor of $i\xi$ from the spatial derivative and we will take fewer time derivatives: $\frac{2s-1}{4}$ instead of $\frac{2s+1}{4}$. 

To estimate $\widetilde{\I}$, notice that on the region of integration $\tau \approx \sqrt{\xi^2 + \xi^4}$, so the additional $|\xi|$ factor is equivalent to $|\tau|^{1/2}$. This brings us exactly back to the situation addressed above for $\I$.

To estimate $\widetilde{\III}$, when $|\xi| \lesssim 1$, the bounds are identical to those for $\III$. When $|\xi| \gg1$, we change variables as we did for $\III$. The additional factor of $\xi$ is equivalent to a factor of $|z|^{1/2}$, which exactly replaces the lost time derivative, and we are again back to the situation addressed in bounding $\III$. 

Estimating $\widetilde{\II}$ is a bit more complex. If $s \geq \frac12$, we have  
\begin{align*}
&\| \widetilde{\II} \|_{H^\frac{2s-1}{4}_t} \; \lesssim \;\; \Bigl\| \lb \tau \rb^{\frac{2s-1}4} \int \frac{|\xi|}{\lb \tau \mp \sqrt{\xi^2 + \xi^4} \rb} \frac{| \widehat{G}(\xi,\tau)|}{\sqrt{\xi^2 +\xi^4}} \d \xi \Bigr\|_{L^2_\tau} \\
&\lesssim \Bigl\| \int \chi_R(\xi, \tau) \lb |\tau| - \xi^2 \rb^\frac{2s-3}{4} \frac{| \widehat{G}(\xi,\tau)|}{\sqrt{\xi^2 +\xi^4}} \d \xi \Bigr\|_{L^2_\tau} 
+ \Bigl\| \int \frac{|\xi|^{s + 1/2}}{\lb \tau \mp \sqrt{\xi^2 + \xi^4} \rb} \frac{| \widehat{G}(\xi,\tau)|}{\sqrt{\xi^2 +\xi^4}} \d \xi \Bigr\|_{L^2_\tau}.
\end{align*}
Thus when $s \geq \frac12$, we have the bound
\[ \|\mathcal{M}(G)\|_{X^{s,-b}} + \Bigl\| \int \chi_R(\xi,\tau) \lb |\tau| -  \xi^2\rb^\frac{2s-3}{4} | \widehat{\mathcal{M}(G)}(\xi,\tau)| \d \xi \Bigr\|_{L^2_\tau} \] 
as before. For $s \leq \frac12$, go back to 
\[ \Bigl\| \lb \tau \rb^{\frac{2s-1}4} \int \frac{|\xi|}{\lb \tau \mp \sqrt{\xi^2 + \xi^4} \rb} \frac{| \widehat{G}(\xi,\tau)|}{\sqrt{\xi^2 +\xi^4}} \d \xi \Bigr\|_{L^2_\tau}. \]
On the region where $|\tau| \ll \xi^2$ and $|\xi | \gtrsim 1$, we obtain the bound
\[ \Bigl\| \lb \tau \rb^{\frac{2s-1}4} \int \chi_{Q}(\xi,\tau) \frac{1}{\lb \xi \rb} \frac{| \widehat{G}(\xi,\tau)|}{\sqrt{\xi^2 + \xi^4}} \d \xi \Bigr\|_{L^2_\tau},\]
where $Q = \{|\tau| \ll |\xi|^2\} \cap \{ | \xi| \gtrsim 1\}$. 

On the region where $|\tau| \gtrsim |\xi|^2$ or $|\xi | \lesssim 1$, we have 
\begin{align*}
\| \widetilde{\II} \|_{H^\frac{2s-1}{4}_t} &\lesssim \Bigl\| \lb \tau \rb^{\frac{2s-1}4} \int \frac{|\xi|}{\lb \tau \mp \sqrt{\xi^2 + \xi^4} \rb} \frac{| \widehat{G}(\xi,\tau)|}{\sqrt{\xi^2 +\xi^4}} \d \xi \Bigr\|_{L^2_\tau} \\
&\lesssim \Bigl\| \int \frac{|\xi|^{s + 1/2}}{\lb \tau \mp \sqrt{\xi^2 + \xi^4} \rb} \frac{| \widehat{G}(\xi,\tau)|}{\sqrt{\xi^2 +\xi^4}} \d \xi \Bigr\|_{L^2_\tau},
\end{align*}
which can be bounded by $\|\mathcal{M}G\|_{X^{s,-b}}$ as we have already seen. This completes the proof.

\subsection{Proof of Lemma \ref{bilinear0}}\label{correction_prf}

We want to show that 
\begin{equation*}
 \left\| \lb \tau \rb^{\frac{2(s+a)-1}{4}} \int_{\substack{|\xi| \gtrsim 1 \\ |\tau| \ll \xi^2}}  \frac{\xi^2 \; \widehat{uv}(\xi,\tau)}{ \lb \xi \rb \sqrt{\xi^2 + \xi^4}} \d \xi \right\|_{L^2_\tau} \lesssim \| u \|_{X^{s,b}} \|v\|_{X^{s,b}}
\end{equation*}
for $\frac12 - b > 0$ sufficiently small. 

Writing the Fourier transform of $uv$ as a convolution, we have 
\[ \widehat{uv}(\xi, \tau) = \iint \widehat{u}(\xi_1, \tau_1) \widehat{v}(\xi -\xi_1, \tau - \tau_1) \d \xi_1 \d \tau_1. \]
Let $f(\xi, \tau) = \widehat{u}(\xi,\tau) \lb \xi \rb^s \lb |\tau| - \xi^2 \rb^b $ and $g(\xi,\tau) = \widehat{v}(\xi,\tau) \lb \xi \rb^s \lb |\tau| - \xi^2 \rb^b$. Using this and dropping the factor $\xi^2/\sqrt{\xi^2 + \xi^4}$, the desired bound becomes
\begin{equation} \label{eq:estimate}
\begin{split}
 \left\| \lb \tau \rb^{\frac{2(s+a)-1}{4}} \iiint_{\substack{|\xi| \gtrsim 1 \\ |\tau| \ll \xi^2}} \frac{f(\xi_1, \tau_1) g(\xi - \xi_1, \tau - \tau_1 ) \lb \xi_1 \rb^{-s} \lb \xi- \xi_1 \rb^{-s}}{ \lb \xi \rb \lb |\tau_1| - \xi_1^2 \rb^b \lb |\tau-\tau_1| - (\xi - \xi_1)^2 \rb^b} \d \xi \d \xi_1 \d \tau_1 \right\|_{L^2_\tau} \qquad \\
 \lesssim \| f\|_{L^2_\xi L^2_\tau} \| g\|_{L^2_\xi L^2_\tau}. 
 \end{split}
\end{equation}

\noindent \textbf{Case 1.} $\operatorname{sgn} (\tau_1) = \operatorname{sgn}(\tau - \tau_1)$. \\
Using the Cauchy-Schwarz inequality in the $\xi$-$\xi_1$-$\tau_1$ integral of \eqref{eq:estimate}, then Cauchy-Schwarz in $\tau$, and finally Young's inequality, we obtain the bounds
\begin{align*} 
& \left\| \| M \|_{L^2_{\xi,\xi_1,\tau_1}\bigl(|\xi| \gtrsim 1, |\tau| \ll \xi^2\bigr)} \| f(\xi_1, \tau_1) g(\xi - \xi_1, \tau - \tau_1) \|_{L^2_{\xi,\xi_1, \tau_1}}\right\|_{L^2_\tau} \\
&\lesssim \left( \sup_\tau \iiint_{\substack{|\xi| \gtrsim 1 \\ |\tau| \ll \xi^2}} M^2 \d \xi \d \xi_1 \d \tau_1 \right)^{1/2} \| f^2 \ast g^2 \|_{L^1_{\xi,\tau}}^{1/2} \\
&\lesssim \left( \sup_\tau \iiint_{\substack{|\xi| \gtrsim 1 \\ |\tau| \ll \xi^2}} M^2 \d \xi \d \xi_1 \d \tau_1 \right)^{1/2} \|f\|_{L^2_{\xi,\tau}} \| g\|_{L^2_{\xi,\tau}}. 
\end{align*}
where 
\[ M = M(\xi_1, \xi, \tau, \tau_1) = \frac{\lb \tau \rb^{\frac{2(s+a) -1}{4}}}{\lb \xi \rb \lb \xi_1 \rb^{s} \lb \xi - \xi_1 \rb^{s} \lb |\tau_1| - \xi_1^2 \rb^b \lb |\tau - \tau_1| - (\xi-\xi_1)^2 \rb^b}. \]
Thus, it suffices to show that the supremum above is finite. Using Lemma \ref{calc_est} in the $\tau_1$ integral, along with the assumption that $\tau_1$ and $\tau - \tau_1$ have the same sign, we arrive at 
\begin{align*} 
 \sup_\tau \iint_{\substack{|\xi| \gtrsim 1 \\ |\tau| \ll \xi^2}} \frac{ \lb \tau \rb^{\frac{2(s+a) -1}{2}} \lb \xi \rb^{-2s} \lb \xi - \xi_1 \rb^{-2s} }{ \lb \xi \rb^2 \lb \tau \pm ( \xi_1 ^2 + (\xi - \xi_1)^2) \rb^{4b-1}} \d \xi \d \xi_1.
\end{align*}
Since $|\tau| \ll \xi^2$ and $\xi_1^2 + (\xi - \xi_1)^2 \gtrsim \max\{ \xi_1^2, \xi^2\}$, we have 
\[\lb \tau \pm ( \xi_1 ^2 + (\xi - \xi_1)^2) \rb^{4b-1} \approx \lb \max\{ |\xi|, |\xi_1|\}  \rb^{8b-2}.\] 
Using this, and dropping the $\lb \tau \rb$ term, it suffices to bound
\[ \iint \frac{\lb \xi_1\rb^{-2s} \lb \xi - \xi_1 \rb^{-2s}}{\lb \xi \rb^2 \lb \max\{ |\xi|, |\xi_1|\} \rb ^{8b-2}} \d \xi \d \xi_1. \]
If $|\xi| \gtrsim |\xi_1|$, this is bounded by 
\begin{align*}
\iint_{|\xi_1| \lesssim |\xi|} \lb \xi_1 \rb^{-2s} \lb \xi \rb^{-8b} \lb \xi - \xi_1 \rb^{-2s} \d \xi \d \xi_1 
\;&\lesssim\; \iint_{|\xi_1| \lesssim |\xi|} \lb \xi_1 \rb^{-2s} \lb \xi \rb^{-8b} \lb \xi \rb^{\max\{-2s, 0\}} \d \xi \d \xi_1 \\
&\lesssim  \int \lb \xi \rb^{1-2s} \lb \xi \rb^{-8b} \lb \xi \rb^{\max\{-2s, 0\}} \d \xi,
\end{align*}
which is finite if $1 - 2s - 8b + \max\{-2s, 0\} < -1$, which holds for $s > \frac12 - 2b$, i.e. $s > -\frac12$ for $\frac12 - b > 0$ sufficiently small. 
If $|\xi| \ll |\xi_1|$, we need to bound
\[ \iint \lb \xi_1 \rb^{2-8b -4s} \lb \xi \rb^{-2} \d \xi \d \xi_1 \; \lesssim \; \int \lb \xi_1 \rb^{2-8b-4s} \d \xi_1,\]
which is finite if $s > -\frac14$ and $\frac12 - b > 0$ is sufficiently small. \\

\noindent \textbf{Case 2.} $\operatorname{sgn} (\tau_1) \neq \operatorname{sgn}(\tau - \tau_1)$ and $| \xi_1 | \lesssim |\xi|$.  \\
Using Cauchy-Schwarz and Young's inequalities just as in Case 1 and dropping the $\lb \tau \rb$ term, it suffices to show that 
\[ \sup_\tau \iint_{\substack{|\xi| \gtrsim 1 \\ |\xi| \gtrsim |\xi_1| \\ |\tau| \ll \xi^2}}\frac{\lb \xi_1 \rb^{-2s} \lb \xi - \xi_1 \rb^{-2s}}{\lb \xi \rb^{2} \lb \tau \pm \xi (\xi - 2 \xi_1)\rb^{2b}} \d \xi_1 \d\xi \]
is finite. Using the change of variables $z = \xi(\xi - 2\xi_1)$ in the $\xi_1$ integral, we arrive at 
\[ \sup_\tau \iint_{\substack{|z| \lesssim \xi^2 \\ |\tau| \ll \xi^2}} \frac{\lb \xi - z/\xi \rb ^{-2s} \lb \xi + z/\xi \rb^{-2s}}{\lb \xi \rb^{3} \lb \tau \pm z \rb^{2b}} \d z \d \xi. \]
Notice that since $|z| \lesssim \xi^2$, we have $|\xi \pm z/\xi| \lesssim |\xi|$. This yields 
\[ \sup_\tau \iint_{\substack{|z| \lesssim \xi^2 \\ |\tau| \ll \xi^2}} \frac{\lb \xi \rb ^{\max\{-4s,0\}}}{\lb \xi \rb^{3} \lb \tau \pm z \rb^{2b}} \d z \d \xi \; \lesssim \; \int \lb \xi \rb^{\max\{-4s, 0\}} \lb \xi \rb^{-3} \lb \xi \rb^{2(1-2b)} \d \xi, \]
which is finite for $s > -\frac12$ and $\frac12 - b >0$ sufficiently small. \\

\noindent \textbf{Case 3.} $\operatorname{sgn} (\tau_1) \neq \operatorname{sgn}(\tau - \tau_1)$ and $| \xi| \ll |\xi_1|$.  \\
By duality, to establish \eqref{eq:estimate}, it suffices to show that 
\[ \left| \iiiint_{\substack{|\tau| \ll \xi^2 \\ 1 \lesssim |\xi| \ll |\xi_1|}} M \; f(\xi_1, \tau_1) g(\xi-\xi_1, \tau - \tau_1) \; \phi(\tau) \d \xi \d \tau \d \xi_1 \d \tau_1 \right| \lesssim \| \phi \|_{L^2_\tau} \|f\|_{L^2_{\xi,\tau}} \|g\|_{L^2_{\xi,\tau}}, \]
where $M = M(\xi,\xi_1, \tau,\tau_1)$ is defined as in Case 1. Using Cauchy-Schwarz in the $\xi_1$-$\tau_1$ integrals, it suffices to show that 
\[ \left\| \iint_{\substack{|\tau| \ll \xi^2 \\ 1 \lesssim |\xi| \ll |\xi_1|}} M \; g(\xi-\xi_1, \tau - \tau_1) \; \phi(\tau) \d \xi \d \tau  \right\|_{L^2_{\xi_1, \tau_1}} \lesssim \| \phi \|_{L^2_\tau} \|g\|_{L^2_{\xi,\tau}}.\]
Using Cauchy-Schwarz in $\xi$-$\tau$ and then Young's inequality, the left-hand side of this quantity is bounded by 
\begin{align*}
 \left\| \| M \lb \xi \rb^{\frac12+} \|_{L^2_{\xi,\tau}\bigl(|\tau| \ll \xi^2, 1 \lesssim |\xi| \ll |\xi_1|\bigr)} \| g(\xi - \xi_1, \tau - \tau_1) \phi(\tau) \lb \xi \rb^{-\frac12 -} \|_{L^2_{\xi,\tau}} \right\|_{L^2_{\xi_1, \tau_1}} \\
 \lesssim \left( \sup_{\xi_1, \tau_1} \iint_{\substack{|\tau| \ll \xi^2 \\  1 \lesssim |\xi| \ll |\xi_1|}} M^2 \lb \xi \rb^{1+} \d \xi \d \tau \right)^{1/2} \| g\|_{L^2_{\xi, \tau}} \| \lb \xi \rb^{-\frac12 -} \|_{L^2_\xi} \| \phi\|_{L^2_\tau}. 
\end{align*}
Thus it suffices to show the following supremum is finite: 
\begin{align*}  &\sup_{\xi_1} \iint_{\substack{|\tau| \ll \xi^2 \\ 1 \lesssim |\xi| \ll |\xi_1|}} \frac{ \lb \tau \rb^{\frac{2(s+a)-1}{2}} \lb \xi_1 \rb^{-2s} \lb \xi - \xi_1 \rb^{-2s}}{\lb \xi \rb^{1-} \lb \tau \pm \xi(\xi - 2 \xi_1) \rb^{2b}} \d \xi \d \tau \\
\lesssim &\sup_{\xi_1}\;  \lb \xi_1 \rb^{-4s} \iint_{\substack{|\tau| \ll \xi^2 \\ |\xi| \ll |\xi_1|}} \frac{ \lb \tau \rb^{\frac{2(s+a)-1}{2}}}{\lb \xi \rb^{1-} \lb \xi \rb^{2b} \lb  \xi_1 \rb^{2b}} \d \xi \d \tau \\
\lesssim &\sup_{\xi_1} \;  \lb \xi_1 \rb^{-2b -4s} \int_{|\xi| \ll |\xi_1|} \frac{\lb \xi \rb^{2(s+a) + 1}}{\lb \xi \rb^{1-} \lb \xi \rb^{2b}} \d \xi \\
\lesssim &\sup_{\xi_1} \;  \lb \xi_1 \rb^{-2b -4s} \lb \xi_1 \rb^{\max\{0, 2(s+a) - 2b  + 1 +\}}.
\end{align*}
If the maximum in the last line is zero, we have a finite bound for $s > - \frac14$ if $\frac12 - b >0$ is sufficiently small. Otherwise we require $ -4b - 4s + 2(s+a) + 1 < 0$, which holds for $a < s + \frac12$ as long as $\frac12 - b > 0$ is sufficiently small.

\subsection{Proof of Lemma \ref{bilinear}} \label{bilinear_prf} 

Recall that we want to show that for $ \frac12 < s + a \leq \frac52$ and $a < \min\{1, s + \frac12\}$, we have 
\[ \Bigl\| \int \chi_R(\xi,\tau) \lb |\tau| -  \xi^2 \rb^\frac{2(s+a)-3}{4} \frac{\xi^2}{\sqrt{\xi^2 + \xi^4}} |\widehat{uv}(\xi,\tau)| \d \xi \Bigr\|_{L^2_\tau} \lesssim \|u\|_{X^{s,b}}\|v\|_{X^{s,b}}.  \]

Writing the Fourier transform as a convolution and canceling the $\xi^2/\sqrt{\xi^2 + \xi^4}$ factor, we need to bound
\begin{align*}
&\Bigl\| \iiint \chi_R(\xi,\tau) \lb |\tau| - \xi^2 \rb^\frac{2(s+a)-3}{4} |\widehat{u}(\xi_1, \tau_1)|| \widehat{v}(\xi -\xi_1, \tau - \tau_1)| \d \xi_1 \d \tau_1 \d \xi  \Bigr\|_{L^2_\tau} \\
&= \Bigl\| \iiint \frac{ \chi_R(\xi,\tau) \;  \lb |\tau| - \xi^2  \rb^\frac{2(s+a)-3}{4}|f(\xi_1, \tau_1)|| g(\xi -\xi_1, \tau - \tau_1)|\; \d \xi_1 \d \tau_1 \d \xi }{\lb \xi_1 \rb^{s} \lb \xi - \xi_1 \rb^{s} \lb |\tau_1| - \xi_1^2  \rb ^b \lb | \tau -\tau_1| - (\xi - \xi_1)^2  \rb^b }    \Bigr\|_{L^2_\tau}
\end{align*}
where 
\[ f(\xi, \tau) = \widehat{u}(\xi, \tau) \lb \xi \rb^s  \lb |\tau| - \xi^2  \rb^{b} \qquad g(\xi, \tau) = \widehat{v}(\xi, \tau) \lb \xi \rb^s  \lb |\tau| - \xi^2 \rb^{b}. \] 
Using the Cauchy-Schwarz inequality in the $\xi_1$-$\tau_1$-$\xi$ integral, followed by Young's inequality as in the proof of Lemma \ref{bilinear0}, it suffices to show that 
\[ \sup_\tau \; \iiint \chi_R(\xi,\tau) \frac{ \lb |\tau|- \xi^2  \rb^\frac{2(s+a)-3}{2}\lb \xi_1 \rb^{-2s} \lb \xi - \xi_1 \rb^{-2s}}{ \lb |\tau_1| - \xi^2_1  \rb ^{2b} \lb | \tau -\tau_1| - (\xi - \xi_1)^2  \rb^{2b} }  \d \xi_1 \d \tau_1 \d \xi \]
is finite. If $\tau_1$ and $\tau - \tau_1$ have the same sign, we apply Lemma \ref{calc_est} in the $\tau_1$ integral and obtain the bound
\begin{equation} \label{eq:q1}
\sup_\tau \; \iint \chi_R(\xi,\tau) \frac{ \lb |\tau| - \xi^2 \rb^\frac{2(s+a)-3}{2}\lb \xi_1 \rb^{-2s} \lb \xi - \xi_1 \rb^{-2s}}{ \lb |\tau| - \xi^2 + 2\xi_1(\xi-\xi_1) \rb ^{2b} }  \d \xi_1  \d \xi. 
\end{equation}
If $\tau_1$ and $\tau - \tau_1$ have different signs, it's bounded by
\[
\sup_\tau \; \iint \chi_R(\xi,\tau) \frac{ \lb \lambda \tau - \xi^2 \rb^\frac{2(s+a)-3}{2}\lb \xi_1 \rb^{-2s} \lb \xi - \xi_1 \rb^{-2s}}{ \lb \lambda \tau - \xi^2 + 2\xi\xi_1) \rb ^{2b} }  \d \xi_1  \d \xi, 
\]
where $\lambda = \operatorname{sgn}(\tau - \tau_1) = \pm 1$. Here we have taken advantage of the fact that we're confined to the set $R$ to conclude that $\lb |\tau| - \xi^2 \rb \approx \lb \lambda \tau - \xi^2 \rb $. Changing variables in the $\xi_1$ integral by  $\xi_1 \mapsto \xi - \xi_1$, and dropping the $\lambda$, we obtain
\begin{equation}\label{eq:q2}
\sup_\tau \; \iint \chi_R(\xi,\tau) \frac{ \lb \tau - \xi^2 \rb^\frac{2(s+a)-3}{2}\lb \xi - \xi_1 \rb^{-2s} \lb  \xi_1 \rb^{-2s}}{ \lb \tau - \xi^2 + 2 \xi(\xi -\xi_1) \rb ^{2b} }  \d \xi_1  \d \xi. 
\end{equation}

When $\frac32 \leq s+a < \frac52$, we use the inequalities
\begin{align*}
\lb |\tau| - \xi^2  \rb &\lesssim \lb |\tau| - \xi^2 + 2 \xi_1(\xi - \xi_1) \rb \lb \xi_1 \rb \lb \xi - \xi_1 \rb \quad \text{and} \\
\lb \tau - \xi^2  \rb &\lesssim \lb \tau - \xi^2 + 2 \xi(\xi - \xi_1)   \rb \lb \xi   \rb \lb \xi - \xi_1 \rb
\end{align*}
in \eqref{eq:q1} and \eqref{eq:q2} respectively. They yield the following bounds for \eqref{eq:q1} and \eqref{eq:q2}:  
\begin{align*}
&\sup_\tau \; \iint_{R} \frac{1}{\lb \xi_1 \rb^{s + \frac32 - a} \lb \xi - \xi_1 \rb ^{s + \frac32 - a}} \d \xi_1 \d \xi &\text{ for \eqref{eq:q1} and} \\
&\sup_\tau \;   \iint_{R} \frac{\lb \xi \rb^{s + a - \frac32}}{\lb \xi_1 \rb^{2s} \lb \xi - \xi_1 \rb ^{s + \frac32 - a}} \d \xi_1 \d \xi &\text{ for \eqref{eq:q2}. }
\end{align*}
Using Lemma \ref{calc_est}, we see that former is finite as long as $s + \frac32 -a > 1$, which holds when $a < s + \frac12$. For the latter, we use Lemma \ref{calc_est} in the $\xi_1$ integral, using the assumption that $a < s + \frac12$, to obtain $\int \lb \xi_1 \rb^{s+ a - \frac32 - (s + \frac32 - a)} \d \xi_1$, which is convergent for $a < 1$. 

When $\frac12 < s +a < \frac32$, we use the inequality $\lb \tau - a \rb \lb \tau - b \rb \gtrsim \lb a-b \rb$ to obtain
\begin{align*}
&\iint \frac{1}{\lb \xi_1 \rb^{2s} \lb \xi - \xi_1 \rb ^{2s}\lb \xi_1(\xi - \xi_1)\rb^{\frac32 - s -a}} \d \xi_1 \d \xi &\text{ from \eqref{eq:q1} and} \\
&\iint \frac{1}{\lb \xi_1 \rb^{2s} \lb \xi - \xi_1 \rb ^{2s } \lb \xi(\xi_1 - \xi) \rb^{\frac32 - s -a}} \d \xi_1 \d \xi &\text{ from \eqref{eq:q2}. }
\end{align*}
In the nonresonant cases, i.e. when $|\xi_1|, |\xi - \xi_1| \gtrsim 1$ for the first equation and when $|\xi|, |\xi - \xi_1| \gtrsim 1$ for the second equation, we have  $\lb \xi_1(\xi - \xi_1)\rb \approx \lb \xi_1 \rb \lb \xi - \xi_1 \rb$ and $\lb \xi(\xi_1 - \xi) \rb \approx \lb \xi \rb \lb \xi_1 - \xi \rb$ respectively. Thus we have convergence if $a < s + \frac12$ for the first equation. In the second equation, we use Lemma \ref{calc_est} to the estimate the $\xi_1$ integral. This yields a bound of $\int \lb \xi \rb^{\frac32 + s -a } \d \xi$, which is convergent if  $a < s+ \frac12$ . 

The resonances can be treated simply. In the first equation, when $|\xi_1| \lesssim 1$, we have 
\[ \iint \frac{1}{\lb \xi_1 \rb^{2s} \lb \xi - \xi_1 \rb ^{2s}\lb \xi_1(\xi - \xi_1)\rb^{\frac32 - s -a}} \d \xi_1 \d \xi \lesssim \int_{-1}^1 \int \lb \xi - \xi_1 \rb^{-2s} \d \xi \d \xi_1.\]
This converges since $s > \frac12$. The remaining resonances can be handled in exactly the same way -- drop two of the three factors, and integrate, using the fact that we're integrating over a finite interval in one of the dimensions and that $s > \frac12$ to obtain convergence.

\section{Appendix}\label{appendix}

\subsection{Proof of Lemma \ref{explicit_sol}: Explicit Linear Solution Formula}

Denote the Laplace transform of a function $u(t)$ defined on $[0,\infty)$ by 
\[ \wt{u}(\lambda) = \int_0^\infty e^{-\lambda t} u(t) \d t. \]
Taking the Laplace transform in time  of \eqref{eq:zB} yields the equation
\begin{equation*}
 \begin{cases}
 \lambda^2 \wt{v}(x,\lambda) - \wt{v}_{xx}(x,\lambda) + \wt{v}_{xxxx}(x,\lambda) = 0, \\
 \wt{v}(0,\lambda) = \wt{h_1}(\lambda) \quad \wt{v}_x(0, \lambda) = \wt{h_2}(\lambda). 
 \end{cases}
\end{equation*}
The characteristic equation of this is $\lambda^2 - w^2 + w^4 = 0$, 
which has roots satisfying 
\[ w^2 = \frac12 \pm \sqrt{\frac14 - \lambda^2}. \]
Notice that $\sqrt{\frac14 - \lambda^2}$ can be defined analytically on $\mathbb{C} \backslash [-1/2,1/2]$ by 
\[ \left|\frac14 - \lambda^2\right|^{1/2} e^{i(\theta_1 + \theta_2 + \pi)/2}, \]
where $\theta_1 = \arg(\lambda + \frac12)$ and $\theta_2 = \arg(\lambda - \frac12)$. This map sends 
\[ \{\lambda \in \mathbb{C} : \Re \lambda \geq 0, \;\; \lambda \notin [0,1/2]\} \mapsto \{\lambda \in \mathbb{C} : \Im \lambda \geq 0, \;\; \lambda \notin [-1/2,1/2]\}. \]

Let 
\[ a = - \left( \frac12 + \sqrt{\frac14 - \lambda^2} \right)^{1/2} \qquad  b = - \left( \frac12 - \sqrt{ \frac14 - \lambda^2} \right)^{1/2},\]
where the outermost root in $a$ is defined with a branch cut in the bottom half-plane and the outermost root in $b$ is defined with a branch cut in the top half-plane.

Then $a$ and $b$ are analytic for $\lambda$ in the closed right half-plane except for the branch cut $[-1/2, 1/2]$. We also have $\Re a, \Re b \leq 0$ for all $\lambda$ in the closed right half-plane. Since we're interested in solutions which decay at infinity, we only need concern ourselves with these two roots of the characteristic equation. Thus, supressing the $\lambda$ dependence of $a$ and $b$, we have 
\[ \wt{u}(x, \lambda) = \frac{1}{a - b}\left[ \biggl( a \wt{h_1}(\lambda) - \wt{h_2}(\lambda) \biggr) e^{bx} - \biggl( b \wt{h_1}(\lambda) - \wt{h_2}(\lambda) \biggr) e^{ax} \right].\]
By Mellin inversion, we have, for any $c > \frac12$, the equality
\begin{align*}
 v(x,t) &= \frac{1}{2 \pi i} \int_{c - i \infty}^{c + i \infty} \frac{e^{\lambda t}}{a - b}\left[ \biggl( a \wt{h_1}(\lambda) - \wt{h_2}(\lambda) \biggr) e^{bx} - \biggl( b \wt{h_1}(\lambda) - \wt{h_2}(\lambda) \biggr) e^{ax} \right] \d \lambda\\
 &= \frac{1}{2 \pi i} \int_{c - i \infty}^{c + i \infty} \frac{e^{\lambda t}}{a^2 - b^2}(a + b)\left[ \biggl( a \wt{h_1}(\lambda) - \wt{h_2}(\lambda) \biggr) e^{bx} - \biggl( b \wt{h_1}(\lambda) - \wt{h_2}(\lambda) \biggr) e^{ax} \right] \d \lambda.
\end{align*}

\begin{figure}
\begin{tikzpicture}
\draw [help lines,->] (-0.5, 0) -- (2.0, 0);
\draw [help lines,->] (0, -4.5) -- (0, 4.5);
\draw[line width=1pt,   decoration={ markings,
  mark=at position 0.25 with {\arrow[line width=1.2pt]{>}},
  mark=at position 0.40 with {\arrow[line width=1.2pt]{>}},
  mark=at position 0.55 with {\arrow[line width=1.2pt]{>}},
  mark=at position 0.65 with {\arrow[line width=1.2pt]{>}},
  mark=at position 0.689 with {\arrow[line width=1.2pt]{>}},
  mark=at position 0.73 with {\arrow[line width=1.2pt]{>}},
  mark=at position 0.87 with {\arrow[line width=1.2pt]{>}},
  mark=at position 0.97 with {\arrow[line width=1.2pt]{>}}},
  postaction={decorate}]
  (1.4, -4.0) -- (1.4, 4.0) -- (0.0, 4.0) -- (0.0, 0.2/2) -- (0.8268, 0.2/2) arc (150:-150:0.2) -- (0.8268, -0.2/2) -- (0.0, - 0.2/2) -- (0.0, - 4.0) -- cycle; 

\node at (1.8, -0.2){$ \mathbb{R}$};
\node at (-0.3, 4.3){$i\mathbb{R}$};

\end{tikzpicture}\caption{The contour of integration} \label{contour}
\end{figure}
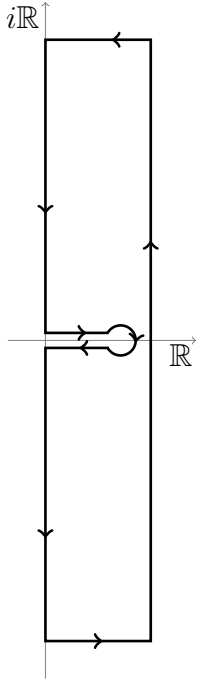

We can write this as an integral along the imaginary axis plus integrals along a keyhole contour about the branch cut and integrals along $s \pm iR$ for $s \in [0 ,c]$ with $R \to \infty$, as shown in Figure \ref{contour}.
The loop of radius $\epsilon$ about the singularity at $\lambda = 1/2$ can be disregarded since the integrand is at most order $1/(a^2 - b^2) \approx |\lambda - 1/2|^{-1/2} \approx \epsilon^{-1/2}$ there, while the length of the contour is order $\epsilon$. 

The integration along the lines $s \pm i\epsilon$ for $s \in [0, 1/2-\epsilon]$ vanishes in the limit as $\epsilon \to 0$ -- the integrals along the two lines cancel one another. This happens because $a(\bar{\lambda}) = \overline{a(\lambda)}$ and $b(\bar{\lambda}) = \overline{b(\lambda)}$. Thus integration over the two lines $s \pm i\epsilon$ for $s \in [0, 1/2-\epsilon]$ is equal to twice the imaginary part of the integral over one of the lines. But the imaginary part of the integrand vanishes as $\epsilon \to 0$. 

The decay of the integrals along $s \pm iR$ for $s \in [0 ,c]$ as $R \to \infty$ is justified as follows. By integration by parts, we have the bound
\[ \left|\wt{h_i}(s \pm iR)\right| \lesssim R^{-1}\Bigl( \| h_i\|_{L^\infty} + \|h_i\|_{L^1} + \|h_i'\|_{L^1} \Bigr) \]
for $s \in [0,c]$. We also have
\[ |a|, |b| \lesssim R^{1/2} \quad \text{and} \quad |a^2 - b^2| \approx R \]
for $\lambda = s \pm iR$ with $R$ large. Thus, on these segments the integrand is order at most $R^{-1}$. Since the intervals are of finite length, we obtain decay as $R \to \infty$. 

Thus we change the contour of integration to the imaginary axis, and arrive at 
\begin{align*}
 v(x,t) &= \frac{1}{2 \pi i} \int_{-i \infty}^{ i \infty} \frac{e^{\lambda t}}{a - b}\left[ \biggl( a \wt{h_1}(\lambda) - \wt{h_2}(\lambda) \biggr) e^{bx} - \biggl( b \wt{h_1}(\lambda) - \wt{h_2}(\lambda) \biggr) e^{ax} \right] \d \lambda\\
 &= \frac{1}{2 \pi i} \int_{ - i \infty}^{ i \infty} \frac{e^{\lambda t}}{a^2 - b^2}(a + b)\left[ \biggl( a \wt{h_1}(\lambda) - \wt{h_2}(\lambda) \biggr) e^{bx} - \biggl( b \wt{h_1}(\lambda) - \wt{h_2}(\lambda) \biggr) e^{ax} \right] \d \lambda \\
&= 2 \Re \frac{1}{2 \pi i} \int_{0}^{ i \infty} \frac{e^{\lambda t}}{a^2 - b^2}(a + b)\left[ \biggl( a \wt{h_1}(\lambda) - \wt{h_2}(\lambda) \biggr) e^{bx} - \biggl( b \wt{h_1}(\lambda) - \wt{h_2}(\lambda) \biggr) e^{ax} \right] \d \lambda.
\end{align*}
Make the change of variables $\lambda = i \mu \sqrt{\mu^2 +1}$. Then $\d \lambda = i\frac{2 \mu ^2 + 1}{\sqrt{\mu^2 + 1}} \d \mu$. On the positive imaginary axis 
\[ a = - i \mu \quad \text{ and } \quad b = - \sqrt{\mu^2 + 1 }, \] 
so $v = \frac{1}\pi \Re(A_0 + B_0 + C_0 +D_0)$, where
\begin{align*}
A_0 &=  -\int_{0}^{ \infty} \frac{e^{i t \mu \sqrt{\mu^2 + 1} -x \sqrt{\mu^2 + 1}}}{\sqrt{1 + \mu^2}}i \mu\biggl(i \mu + \sqrt{1 + \mu^2}\biggr)\; \widehat{h_1}\left( \mu \sqrt{\mu^2 +1}\right) \d \mu \\
B_0 &=  -\int_{0}^{ \infty} \frac{e^{i t \mu \sqrt{\mu^2 + 1} - x\sqrt{\mu^2 + 1}}}{\sqrt{1 + \mu^2}}     \biggl(i \mu + \sqrt{1 + \mu^2}\biggr)\; \widehat{h_2}\left( \mu \sqrt{\mu^2 +1}\right) \d \mu \\
C_0 &=   \int_{0}^{ \infty} e^{i t \mu \sqrt{\mu^2 + 1} - ix\mu}                              \biggl(i \mu + \sqrt{1 + \mu^2}\biggr)            \; \widehat{h_1}\left( \mu \sqrt{\mu^2 +1}\right) \d \mu \\
D_0 &=   \int_{0}^{ \infty} \frac{e^{i t \mu \sqrt{\mu^2 + 1} - ix\mu}}{\sqrt{1 + \mu^2}}     \biggl(i \mu + \sqrt{1 + \mu^2}\biggr)            \; \widehat{h_2}\left( \mu \sqrt{\mu^2 +1}\right) \d \mu. 
\end{align*}

For $x \geq 0$, this is equivalent to $2\pi v(x,t) = -A - B + C + D$. 
Here we used the formula $2\Re z = z + \overline{z}$ to rewrite the real parts of $A_0$, $B_0$, $C_0$, and $D_0$, and added the cut-off function $\rho$ in $A$ and $B$ so that the integrals converge for all $x$. 

\subsection{Calculus Estimates}
The following calculus estimate is important throughout the proofs. See, e.g. \cite{ET3} for a proof. 
\begin{lemma}\label{calc_est}
If $\beta \geq \gamma \geq 0$ and $\beta + \gamma >  1$, then we have
\[ \int \frac{1}{\lb x - a \rb^{\beta} \lb x - b \rb^{\gamma}} \d x \lesssim \lb a - b \rb ^{- \gamma} \varphi_\beta(a-b), \]
where 
\begin{equation*}
 \varphi_\beta(c) = \begin{cases} 1                   & \text{ if } \beta > 1 \\
                                  \log(1 + \lb c \rb) & \text{ if } \beta = 1 \\
                                  \lb c \rb^{1- \beta}& \text{ if } \beta < 1.
                      \end{cases}
\end{equation*}
\end{lemma}

\end{document}